\documentclass[a4paper,12pt]{article}
\usepackage{authblk}

\usepackage{amssymb}
\usepackage{amsmath}

\usepackage{tikz}
\usetikzlibrary{math,decorations.pathreplacing,calligraphy}
\usepackage{diagbox}
\usepackage[width=16.00cm]{geometry}

\def\CX{\cos(k_x \Delta x)}
\def\CY{\cos(k_y \Delta y)}

\def\nn{\nonumber}
\def\refp#1{(\ref{#1})}
\def\showeps{false}

\def\skx{\text{s}(k_x)}
\def\shkx{\text{sh}(k_x)}

\def\shk{\text{sh}(k)}
\def\shfk{\text{sh}'(k)}
\def\shsk{\text{sh}''(k)}

\def\sshk{\text{sh}(k^*)}
\def\sshfk{\text{sh}'(k^*)}

\def\thk{\text{T}(k)}
\def\thfk{\text{T}'(k)}

\def\sthk{\text{T}(k^*)}
\def\sthfk{\text{T}'(k^*)}

\begin{document}

\title{The optimal relaxation parameter for the SOR method applied to the Poisson equation on rectangular grids with different types of boundary conditions}

\author[1]{Hossein Mahmoodi Darian\thanks{hmahmoodi@ut.ac.ir}}

\affil[1]{School of Engineering Science, College of Engineering, University of Tehran}

\date{}

\maketitle

\begin{abstract}
The Successive Over-Relaxation (SOR) method is a useful method for solving the sparse system of linear equations which arises from finite-difference discretization of the Poisson equation. Knowing the optimal value of the relaxation parameter is crucial for fast convergence. In this manuscript, we present the optimal relaxation parameter for the
discretized Poisson equation with mixed and different types of boundary conditions on a rectangular grid with unequal mesh sizes in $x$- and $y$-directions ($\Delta x \neq \Delta y$) which does not addressed in the literature. The central second-order and high-order compact (HOC) schemes are considered for the discretization and the optimal relaxation parameter is obtained for both the point and line implementation of the SOR method. Furthermore, the obtained optimal parameters are verified by numerical results.
\end{abstract}

\subsubsection*{keywords:}
Successive over-relaxation, Optimal relaxation parameter, Poisson equation, Finite difference schemes, High-order compact schemes, Robin boundary conditions

\section{Introduction}
\label{intro}
Numerical solution of Poisson's equation by finite difference schemes appears in many engineering and scientific problems such as heat transfer, fluid mechanics, electrostatics, etc. Iterative methods are used to solve the system of algebraic equations resulting from the discretization of this equation. There are various iterative methods, among which one may mention old methods such as Jacobi and Gauss-Seidel methods \cite{Young1971-va,saad2003iterative}. The convergence rate of the Gauss-Seidel method is twice that of the Jacobi method. However, the convergence rate of both methods is very slow for fine grids. The speed of convergence can be significantly increased by applying the Successive Over-Relaxation (SOR) method  \cite{Young1971-va,saad2003iterative}.
In this method, the difference between the new and old values of each variable is multiplied by a factor and then added to the old value. This factor is called the relaxation parameter. It is very important to know the optimal value of this parameter so that the convergence is achieved in the least number of iterations. It is worth mentioning that there are asymptotically optimal SOR methods \cite{Liu2021,Bai2003,MENG2014707}, however they require more computational cost at each iteration. Therefore, knowing the optimal relaxation parameter is very useful. For the discrete Poisson's equation, the optimal value depends on the grid size and boundary conditions.
However, in the literature the optimal value is mainly presented for square grids with Dirichlet boundary conditions \cite{leveque1988sor,adams1988sor,Tony1989,kuo1990two,YANG2009325} and to the author's knowledge is not presented for rectangular grids where the grid step sizes in $x$- and $y$-directions are not equal ($\Delta x \neq \Delta y$) and other boundary condition types, namely Neumann and Robin boundary conditions. In this article, the optimal relaxation parameter for rectangular grids and different boundary condition types is presented for the use of researchers. In addition to the point SOR method, the optimal parameter is also obtained for the line SOR method. The discretization schemes are the central second-order and high-order compact (HOC) schemes \cite{spotz1995high}. Also, for the latter a more accurate relaxation parameter is obtained in comparison with \cite{leveque1988sor} using the second-order perturbation analysis.

The paper is composed of the following sections. In Section \ref{sec:2}, we give an overview of the discretization schemes and the SOR method. In Section \ref{sec:3}, we obtain the optimal relaxation parameter for all the methods when the boundary conditions are of Dirichlet type. In Sections \ref{sec:4} and \ref{sec:5}, we deal with the Neumann and Robin boundary conditions, respectively. In section \ref{conclusions}, some conclusions are given.

\section{Poisson equation}
\label{sec:2}
The Poisson equation with the Dirichlet boundary conditions in the two-dimensional space is
\begin{align}
	\nabla^2 u \equiv u_{xx}+u_{yy} =&\; f(x,y),	\qquad	(x,y)\in\Omega
	\label{e:laplace}\\
	u(x,y) =&\; g(x,y),	\qquad	(x,y)\in\partial\Omega \nonumber
\end{align}
where $\Omega = [0,1]\times [0,1]$ is the solution domain. Since the convergence speed does not depend on the source term and boundary conditions, we simply set $f(x,y)=g(x,y)=0$. 

\begin{figure}[h]
	\setkeys{Gin}{draft=\showeps}
	\includegraphics[width=0.5\textwidth]{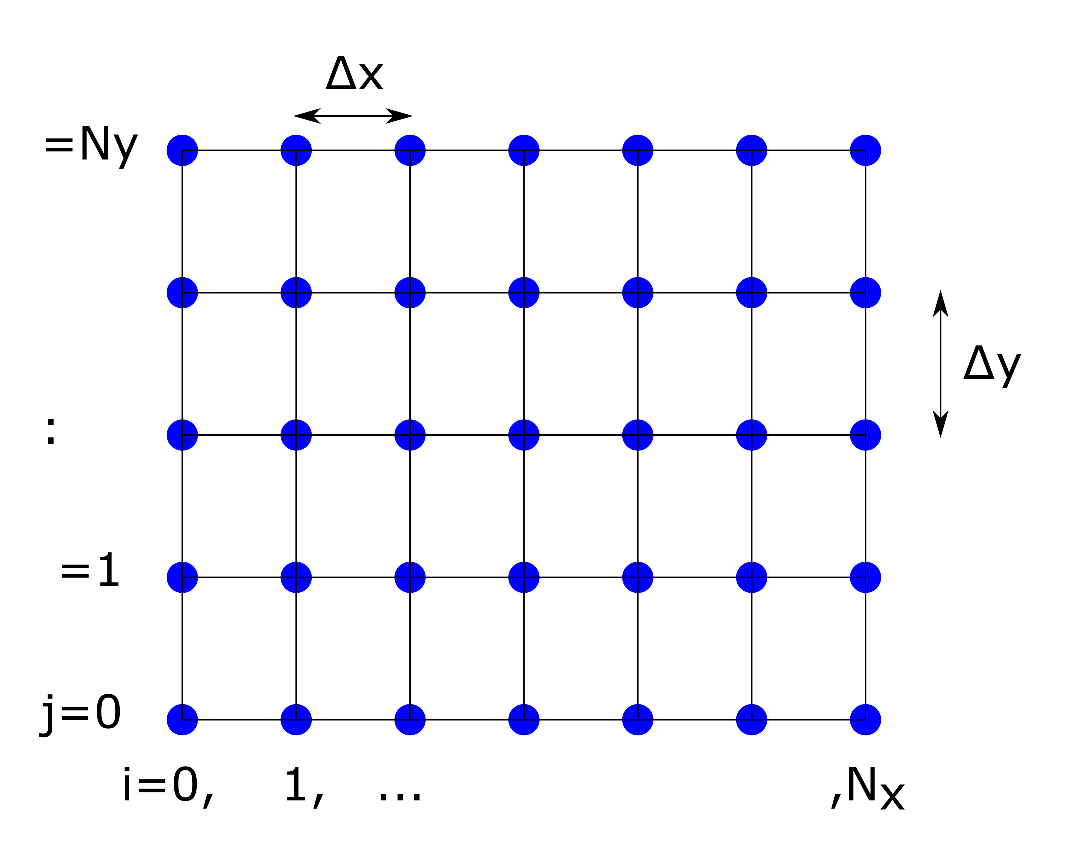}
	\includegraphics[width=0.5\textwidth]{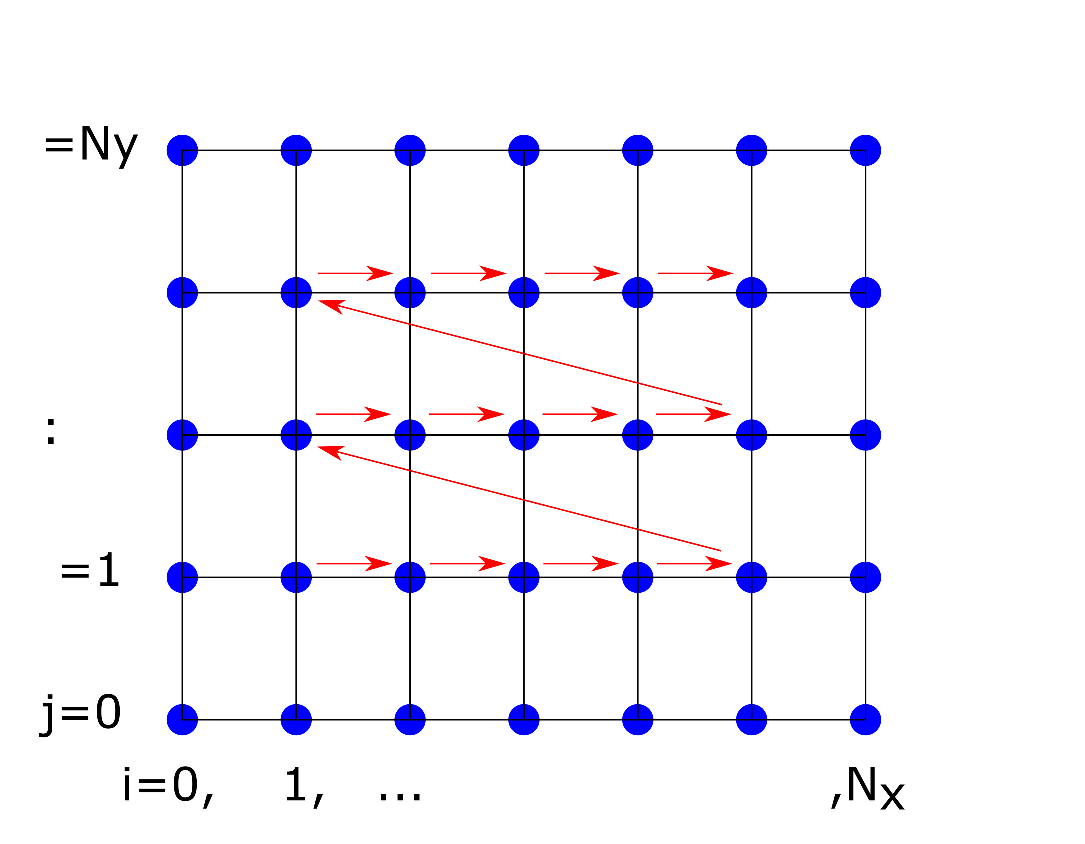}
	\caption{Numerical grid (left) and natural row-wise ordering (Right).}
	\label{fig:grid}       
\end{figure}

A uniform rectangular grid (Fig. \ref{fig:grid}) is used to discretize the domain 
\begin{eqnarray}
	&&x_i = i\Delta x \;,\quad    y_j = j\Delta y \;,\quad   
	u_{i,j} \equiv u(x_i,y_j) \nonumber\\
	&&0\leq i\leq N_x \;, \quad 0\leq j\leq N_y \;,
	\quad \Delta x = \frac{1}{N_x}\;, \quad \Delta y = \frac{1}{N_y}
	\label{e:grid}
\end{eqnarray}
where the number of points in $x$ and $y$ directions are $N_x+1$ and $N_y+1$, respectively.

Discretizing the equation using the central second-order finite difference scheme \cite{hoffmann2000computational} results in
\begin{eqnarray*}
	\frac{u_{i-1,j}-2u_{i,j}+u_{i+1,j}}{\Delta x^2} + \frac{u_{i,j-1}-2u_{i,j}+u_{i,j+1}}{\Delta y^2}= 0
\end{eqnarray*}
where multiplying the both sides by $\Delta x ^ 2$ and defining the grid size ratio as $\beta = \Delta x/\Delta y$, yields
\begin{eqnarray}
	(u_{i-1,j}+u_{i+1,j}) + \beta^2 (u_{i,j-1}+u_{i,j+1}) - 2(1+\beta^2)u_{i,j} = 0 
	\label{e:central}
\end{eqnarray}

Discretization using the HOC scheme \cite{spotz1995high}, gives a 9-point equation
\begin{align}
	&(10-2\beta^2)(u_{i-1,j}+u_{i+1,j}) + (10\beta^2-2) (u_{i,j-1}+u_{i,j+1})  +\nn\\
	&(1+\beta^2)\left[u_{i-1,j-1}+u_{i+1,j-1}+u_{i-1,j+1}+u_{i+1,j+1}- 20 u_{i,j}\right] = 0 
	\label{e:compact}
\end{align}

\subsection{Point SOR method}
\label{sec:pointSOR}
For the Gauss-Seidel method using the natural row-wise ordering (Fig. \ref{fig:grid}), the iteration equations of the second-order and high-order schemes, respectively, are as follows
\begin{align}
	u_{i,j}^{n+1} &= \frac{1}{2(1+\beta^2)}\left((u_{i-1,j}^{n+1}+u_{i+1,j}^{n}) + \beta^2 (u_{i,j-1}^{n+1}+u_{i,j+1}^{n})\right)   
	\label{e:central:gs}\\
	u_{i,j}^{n+1} &=\frac{5-\beta^2}{10(1+\beta^2)}(u_{i-1,j}^{n+1}+u_{i+1,j}^{n}) + \frac{5\beta^2-1}{10(1+\beta^2)} (u_{i,j-1}^{n+1}+u_{i,j+1}^{n})  \nn\\
	&+\frac{1}{20}\left(u_{i-1,j-1}^{n+1}+u_{i+1,j-1}^{n+1}+u_{i-1,j+1}^{n}+u_{i+1,j+1}^{n}\right) 
	\label{e:compact:gs}
\end{align}
where the superscript $n$ is the number of iteration steps.

\begin{figure}[h]
	\setkeys{Gin}{draft=\showeps}
	\begin{center}
		\includegraphics[width=0.6\textwidth]{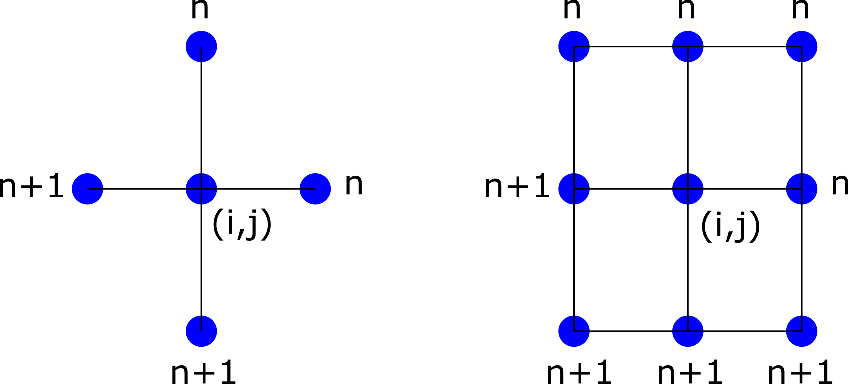}
		\caption{Update state of the neighboring points of $(i,j)$ for the second-order (left) and HOC (right) schemes.}
		\label{fig:order} 
	\end{center}
\end{figure}

Figure \ref{fig:order} shows the update state of the neighboring points of $(i,j)$ when $u_{i,j}$ is updated.

The over-relaxation applies as 
\begin{eqnarray}
	\label{e:relax}
	u_{i,j}^{n+1} = u_{i,j}^{n} + \omega (u_{i,j}^{*}-u_{i,j}^{n}) = (1-\omega)u_{i,j}^{n} + \omega u_{i,j}^{*} 
\end{eqnarray}
where $u_{i,j}^{*}$ is the same as the right-hand side of equations \refp{e:central:gs} or \refp{e:compact:gs} and $\omega$ is the relaxation parameter. The necessary condition for convergence \cite{Axelsson_1994} is $0<\omega<2$. Note that the over-relaxation must be applied immediately after these equations 
not after applying the Gauss-Seidel method to all the points. Otherwise, the optimum convergence rate is not obtained. Due to this, the Gauss-Seidel and over-relaxation steps are combined as a single step:
\begin{align}
	u_{i,j}^{n+1} &= (1-\omega)u_{i,j}^{n}\nn \\&+\frac{\omega}{2(1+\beta^2)}\left((u_{i-1,j}^{n+1}+u_{i+1,j}^{n}) + \beta^2 (u_{i,j-1}^{n+1}+u_{i,j+1}^{n})\right)   
	\label{e:central:sor}
	\\
	u_{i,j}^{n+1} &=(1-\omega)u_{i,j}^{n}\nn
	\\&+\omega\left[\frac{(5-\beta^2)}{10(1+\beta^2)}(u_{i-1,j}^{n+1}+u_{i+1,j}^{n}) + \frac{(5\beta^2-1)}{10(1+\beta^2)} (u_{i,j-1}^{n+1}+u_{i,j+1}^{n})\right.  \nn\\
	&\left.+\frac{1}{20}\left(u_{i-1,j-1}^{n+1}+u_{i+1,j-1}^{n+1}+u_{i-1,j+1}^{n}+u_{i+1,j+1}^{n}\right) \right]
	\label{e:compact:sor}
\end{align}

\subsection{Line SOR method}
\label{sec:lineSOR}
In the line Gauss-Seidel method, all the points of a row ($j=$ const.) are simultaneously updated:
\begin{eqnarray}
	2(1+\beta^2)u_{i,j}^{n+1} -(u_{i-1,j}^{n+1}+u_{i+1,j}^{n+1})  = 
	\beta^2 (u_{i,j-1}^{n+1}+u_{i,j+1}^{n}) 
	\label{e:lgs:central2}
\end{eqnarray}
therefore, for each row a tridiagonal system needs to be solved. Combining \refp{e:lgs:central2} by the relaxation step \refp{e:relax}, gives the line SOR method
\begin{align}
	2(1+\beta^2)u_{i,j}^{n+1} -(u_{i-1,j}^{n+1}+u_{i+1,j}^{n+1})  &= 
	(1 - \omega)\left(2(1+\beta^2)u_{i,j}^{n} -(u_{i-1,j}^{n}+u_{i+1,j}^{n})\right) \nn\\
	&+\omega \beta^2 \left(u_{i,j-1}^{n+1}+u_{i,j+1}^{n}) \right) 
	\label{e:lsor:central2}
\end{align}
and similarly for the HOC scheme, we have
\begin{align}
	&\quad 20(1+\beta^2)u_{i,j}^{n+1} -(10-2\beta^2)(u_{i-1,j}^{n+1}+u_{i+1,j}^{n+1}) \nn\\
	&=(1-\omega)\left[20(1+\beta^2)u_{i,j}^{n}-(10-2\beta^2)(u_{i-1,j}^{n}+u_{i+1,j}^{n})\right] \nn\\
	&+ \omega \left[(10\beta^2-2) (u_{i,j-1}^{n+1}+u_{i,j+1}^{n})\right.  \nn\\
	&+\left.	(1+\beta^2)\left(u_{i-1,j-1}^{n+1}+u_{i+1,j-1}^{n+1}+u_{i-1,j+1}^{n}+u_{i+1,j+1}^{n}\right)\right] 
	\label{e:lsor:hoc}
\end{align}

\section{Optimal relaxation parameter}
\label{sec:3}
In this section we obtain the optimal relaxation factors for the methods given in the previous section. 
The obtained optimal relaxation parameters are verified numerically. Since the source and boundary conditions are zero, the exact solution is also zero ($u_{i,j}=0$).  Therefore, at each point the error is the same as the guessed solution ($u^n_{i,j}$). This allows us to adopt a very small convergence criteria. We take the convergence criteria equal to 
the fourth power of 
the double-precision machine epsilon ($\epsilon$): 
\begin{eqnarray}
	\label{e:conv}
	\|u\|_2 < \epsilon^4 = \left(2^{-52}\right)^4 \simeq 10^{-63} 
\end{eqnarray}
where
\begin{eqnarray*}
	\|u\|_2 = \sqrt{\sum_{i=0}^{N_x}\sum_{j=0}^{N_y} u_{i,j}^2}
\end{eqnarray*}

The initial guess is taken to be one for all the interior points as 
\begin{eqnarray*}
	u_{i,j} = 1 \;, \qquad 1\leq i\leq N_x - 1 \;, \quad 1\leq j\leq N_y - 1
\end{eqnarray*} 
which means there is a discontinuity at the boundaries where $u_{i,j}$ is zero. 

\subsection{second-order scheme}
Assuming the equation \refp{e:central:sor} has a solution of the form
\begin{eqnarray}
	\label{e:eigen}
	u_{i,j}^n = \lambda^n w(x_i,y_j)
\end{eqnarray}
then $\lambda$ is an eigenvalue of the iteration matrix, and the vector with components
$w_{i,j}\equiv w(x_i, y_j) ,\; 1\leq i\leq (N_x - 1),\; 1\leq j\leq (N_y - 1)$ is the corresponding eigenvector. The eigenvalues may be real or complex. The eigenvalue with the largest modulus (spectral radius), determines the convergence rate. Following the method used in \cite{adams1988sor}, putting \refp{e:eigen} into \refp{e:central:sor},
and canceling a common factor of $\lambda^n$, 
gives
\begin{align}	
	\lambda w_{i,j}&=(1-\omega)w_{i,j} \nn\\
	&+ \frac{\omega}{2(1+\beta^2)}\left((\lambda w_{i-1,j}+w_{i+1,j}) + \beta^2 (\lambda w_{i,j-1}+w_{i,j+1})\right)
	\label{e:eigen:sor}
\end{align}

The eigenvector $w_{i,j}$ must be zero on the domain boundary in view of the given boundary conditions.

Note that, the von Neumann method (complex Fourier analysis), due to asymmetry of \refp{e:eigen:sor}, does not satisfy the boundary conditions and therefore is not applicable here \cite{leveque1988sor}.
Therefore, the method of separation of variables is used. Let $w(x_i,y_j) = X(x_i)Y(y_j)\equiv X_iY_j$ and also  $\lambda = \alpha^2$, then substitution into \refp{e:eigen:sor} gives
\begin{eqnarray}
	\label{e:eigen:sor2}
	\frac{\alpha^2+\omega-1}{\omega}=
	\frac{1}{2(1+\beta^2)}
	\left(\frac{\alpha^2 X_{i-1}+X_{i+1}}{X_i}+\beta^2\frac{\alpha^2 Y_{j-1}+Y_{j+1}}{Y_j}\right)
\end{eqnarray}
now we assume 
\begin{equation}
	\label{e:YX}
	Y = \alpha^{\frac{y}{\Delta y}} \sin(k_y y), \qquad X = \alpha^{\frac{x}{\Delta x}} \sin(k_x x)
\end{equation}

The allowable values for $k_x$ and $k_y$ are those which satisfy the boundary conditions, therefore
\begin{align}
	\label{e:kx}
	k_x &= n_x \pi\;,\quad 1 \leq  n_x\leq {N_x-1}  \\
	k_y &= n_y \pi\;, \quad 1 \leq  n_y\leq {N_y-1} 
\end{align}
Furthermore, using triangular identities, we have
\begin{equation*}
	\frac{\alpha^2 Y_{j-1} + Y_{j+1}}{Y_j} = 2\alpha\cos(k_y \Delta y), \qquad 
	\frac{\alpha^2 X_{i-1} + X_{i+1}}{X_i} = 2\alpha\cos(k_x \Delta x)
\end{equation*}
therefore substitution in \refp{e:eigen:sor2} results in
\begin{eqnarray*}
	\frac{\alpha^2 + \omega - 1}{\omega} = \frac{\alpha}{1+\beta^2}
	\left(\cos(\theta_x) +\beta^2\cos(\theta_y)\right)
\end{eqnarray*}
where $\theta_x = k_x\Delta x = \pi n_x/N_x$ and $\theta_y = k_y\Delta y = \pi n_y/N_y$. After rearranging, we obtain
\begin{eqnarray}
	\alpha^2 - r\omega\alpha + \omega - 1 = 0 \; ,\qquad
	r = \frac{\cos(\theta_x) +\beta^2\cos(\theta_y)}{1+\beta^2}	
	\label{e:sor:r}
\end{eqnarray} 
Also, we have $0<\theta_x, \theta_y < \pi$ and $0<r<1$. The roots of \refp{e:sor:r} are
\begin{eqnarray*}
	\alpha_1 = \frac{r\omega+\sqrt{\Delta}}{2}\; , \qquad \alpha_2 = \frac{r\omega-\sqrt{\Delta}}{2}\; , \qquad \Delta = r^2\omega^2-4\omega+4
\end{eqnarray*}
if $\Delta\geq 0$ then the roots are real and $\alpha_1\geq\alpha_2$. Otherwise, the roots are complex conjugate and 
$|\alpha_1|=|\alpha_2|$. For each pair of $(n_x,n_y)$ there are two distinct eigenvalues ($\lambda=\alpha^2$). However, the eigenvalues of the pair $(n_x,n_y)$ are equal to those of the pair $(N_x-n_x,N_y-n_y)$. Therefore, there are $(N_x-1)(N_y-1)$ distinct eigenvalues.

Note that, it can be easily shown that if $r<1$ then $|\alpha_{1}|,|\alpha_{2}|<1$ and if $r>1$ then $|\alpha_1|>1$. It means if $r<1$ for all the allowable values of $k_x$ and $k_y$ then all the eigenvalues are smaller than one in modulus and the method is convergent and if $r>1$ for some $k_x$ and $k_y$, then it is divergent.

When the roots are real, |$\alpha_1$| is a decreasing function of $\omega$ and when the roots are complex, |$\alpha_1$| is an increasing function of $\omega$. Therefore, the value of $\omega$ which minimizes |$\alpha_1$|, occurs when $\Delta =r^2\omega^2-4\omega+4= 0$. This gives
\begin{eqnarray}
	\omega=\frac{2}{1 \mp \sqrt{1-r^2}}
\end{eqnarray}
for the minus sign, we have $\omega>2$ which is outside of the convergence interval of the SOR method. Therefore, only the plus sign must be considered. The largest eigenvalues in modulus corresponds to the largest $\omega$, which occurs for $k_x=k_y=\pi$. Therefore,  
\begin{eqnarray}
	\omega_{opt} = \frac{2}{1 + \sqrt{1-r_{opt}^2}}\;, \qquad
	r_{opt} = \frac{\cos(\pi\Delta x) +\beta^2\cos(\pi\Delta y)}{1+\beta^2}	
	\label{e:opt:2nd}
\end{eqnarray}
also, using $\Delta = 0$, we have
\begin{eqnarray}
	|\alpha_1|^2 = |\alpha_2|^2 = \frac{r_{opt}^2\omega_{opt}^2}{4}=\omega_{opt} - 1
	\label{e:maxrho:2nd}
\end{eqnarray}

\begin{figure}[h]
	\setkeys{Gin}{draft=\showeps}
	\includegraphics[width=0.5\textwidth]{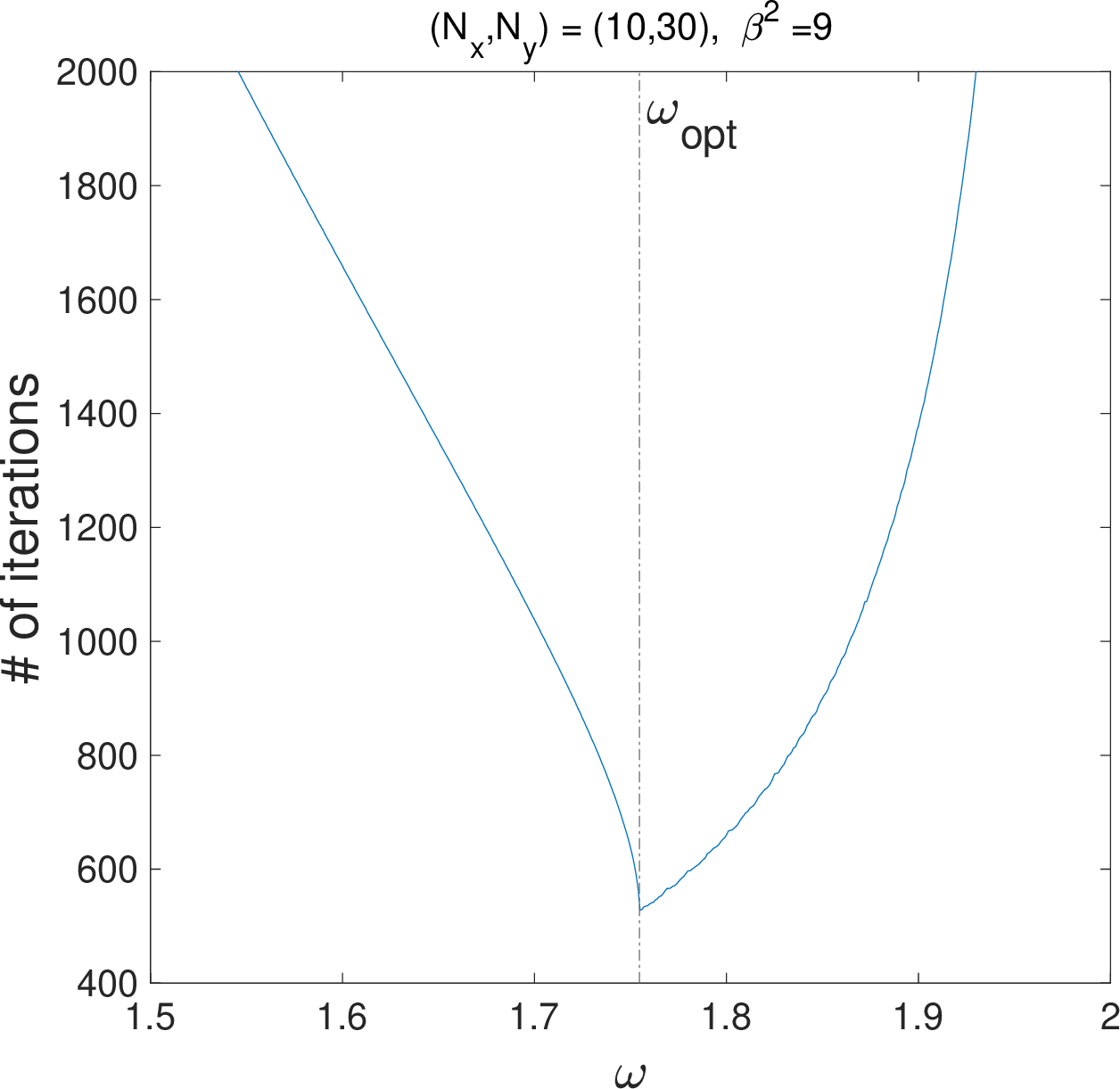}
	\includegraphics[width=0.5\textwidth]{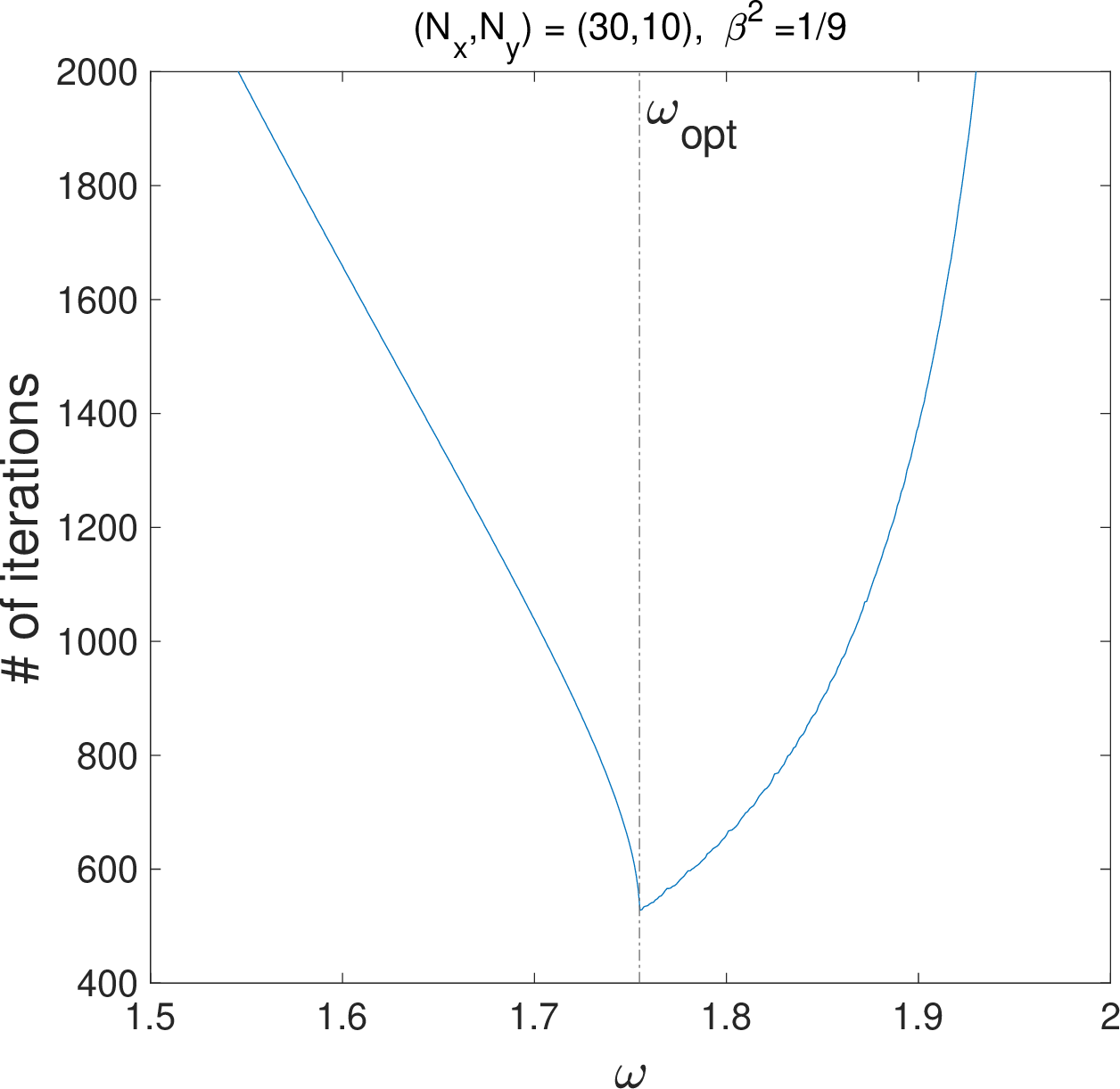}	
	\caption{Number of iterations for the point SOR method of the second-order scheme \refp{e:central:sor} for different values of $\omega$. The dash-dotted line indicates the optimal value \refp{e:opt:2nd}.}
	\label{fig:plot:iter:sor}       
\end{figure}

For numerical verification of the obtained optimal relaxation parameter, we consider two grids with $(N_x,N_y)=(10,30)$ and $(N_x,N_y)=(30,10)$, where the grid size ratios are $\beta = 3$ and $\beta = 1/3$, respectively.
Figure \ref{fig:plot:iter:sor} shows the number of iterations required for satisfying the convergence criterion \refp{e:conv} for different values of $\omega$.  The figure shows the optimal relaxation parameter gives the least number of iterations for both the grids.

For the line SOR method the procedure is similar. The assumed functions for $X(x)$ and $Y(y)$ are as follows: 
\begin{equation*}
	Y = \alpha^{\frac{y}{\Delta y}} \sin(k_y y), \qquad X =  \sin(k_x x)
\end{equation*}
note that the assumed function for $Y(y)$ is identical to the point SOR method, while it is different for $X(x)$. Substitution in \refp{e:lsor:central2} results in a similar equation as \refp{e:sor:r} 
\begin{eqnarray}
	\alpha^2 - r\omega\alpha + \omega - 1 = 0 \; ,\qquad
	r = \frac{\beta^2\cos(\theta_y)}{1+\beta^2-\cos(\theta_x)}	
	\label{e:lsor:r}
\end{eqnarray} 
and therefore
\begin{eqnarray}
	\omega_{opt} = \frac{2}{1 + \sqrt{1-r_{opt}^2}}, \qquad
	r_{opt} = \frac{\beta^2\cos(\pi\Delta y)}{1+\beta^2-\cos(\pi\Delta x)}
	\label{e:opt:lsor}
\end{eqnarray}

\begin{figure}[h]
	\setkeys{Gin}{draft=\showeps}
	\includegraphics[width=0.5\textwidth]{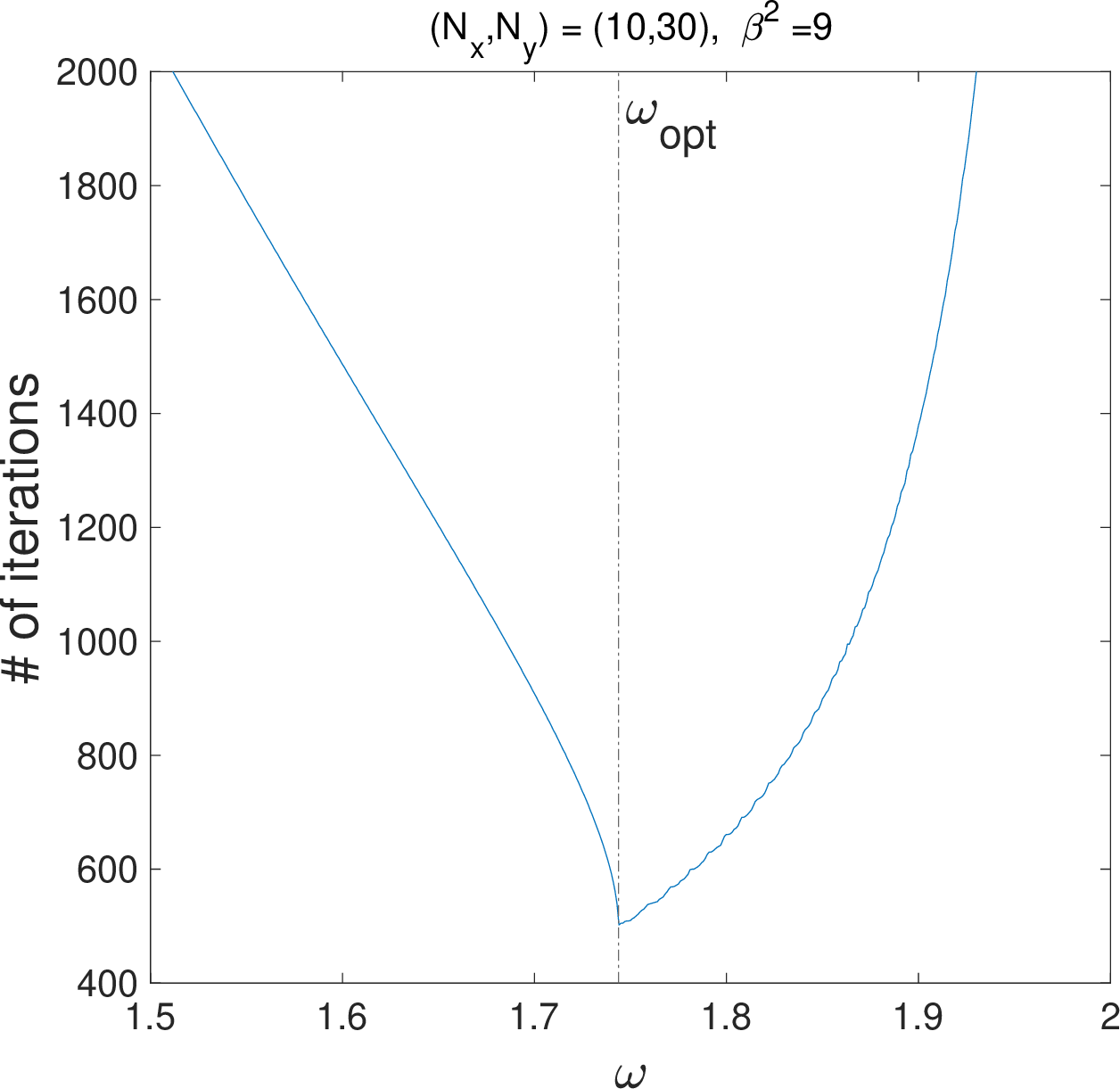}
	\includegraphics[width=0.5\textwidth]{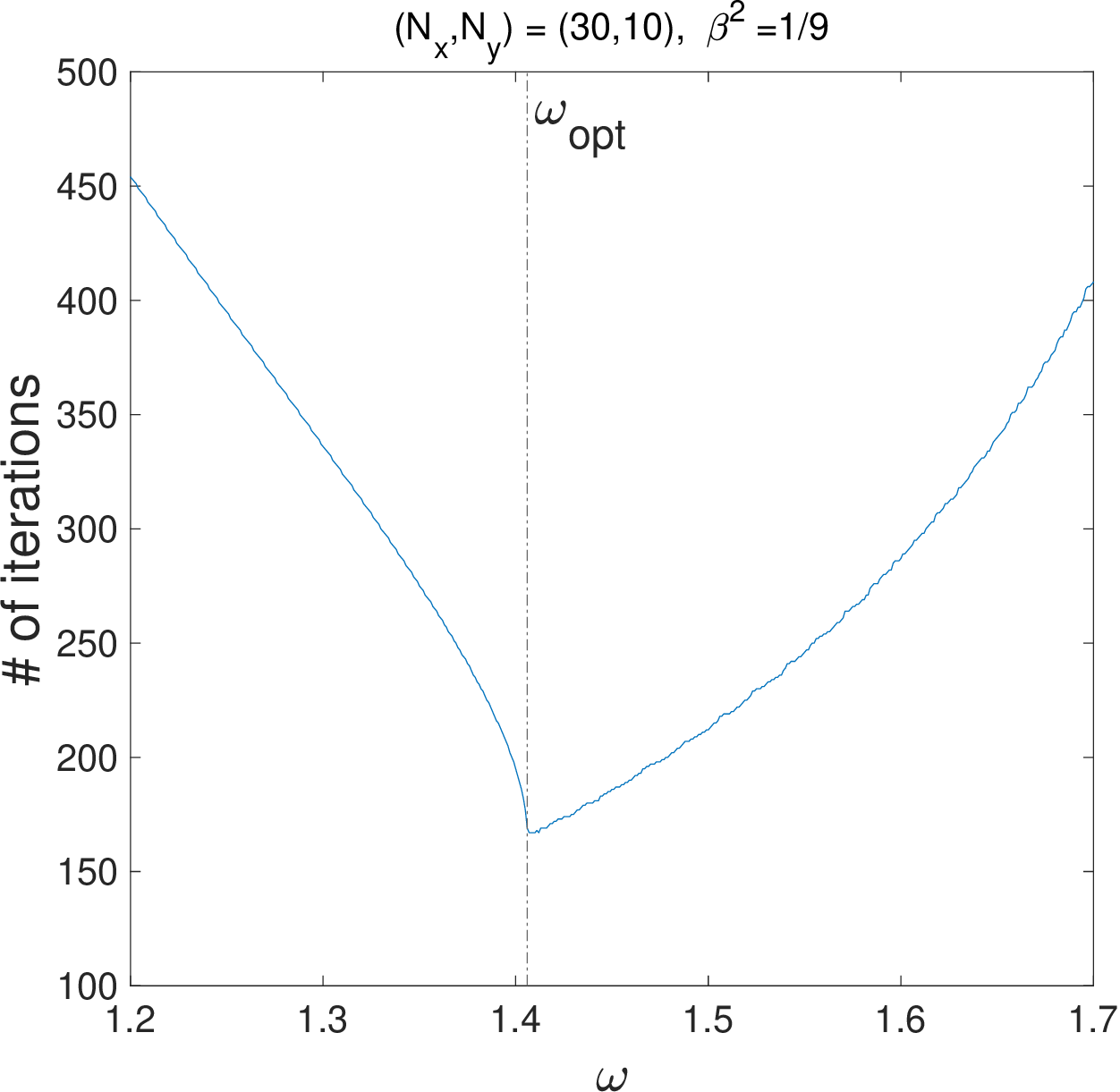}
	\caption{Number of iterations for the line SOR method of the second-order scheme \refp{e:lsor:central2} for different values $\omega$. The dash-dotted line indicates the optimal value \refp{e:opt:lsor}.}
	\label{fig:plot:iter:lsor}       
\end{figure}

Figure \ref{fig:plot:iter:lsor} shows the number of iterations required for satisfying the convergence criterion \refp{e:conv} for different values of $\omega$. The figure shows the optimal relaxation parameter gives the least number of iterations for both the grids. Note that, the number of iterations for the grid with less rows ($N_y=10$) are less then that of the other grid ($N_y=30$). 

\subsection{high-order compact scheme}
Due to the complexity of finding the optimal $\omega$ for the point SOR method of the HOC scheme, we first obtain it for the line SOR method.
The procedure is identical to that of the second-order scheme and a similar equation as \refp{e:lsor:r} is obtained:
\begin{eqnarray}
	\alpha^2 - r\omega\alpha + \omega - 1 = 0 \; ,\quad
	r = \frac{5\beta^2-1+(1+\beta^2)\cos(\theta_x)}{5(1+\beta^2)-(5-\beta^2)\cos(\theta_x)}	\cos(\theta_y)
	\label{e:lsor:hoc:r}
\end{eqnarray} 
therefore
\begin{eqnarray}
	\omega_{opt} = \frac{2}{1 + \sqrt{1-r_{opt}^2}}\;, \;\;
	r_{opt} =\frac{5\beta^2-1+(1+\beta^2)\cos(\pi\Delta x)}{5(1+\beta^2)-(5-\beta^2)\cos(\pi\Delta x)}\cos(\pi\Delta y)
	\label{e:opt:lsor:hoc}
\end{eqnarray}

\begin{figure}[h]
	\setkeys{Gin}{draft=\showeps}
	\includegraphics[width=0.5\textwidth]{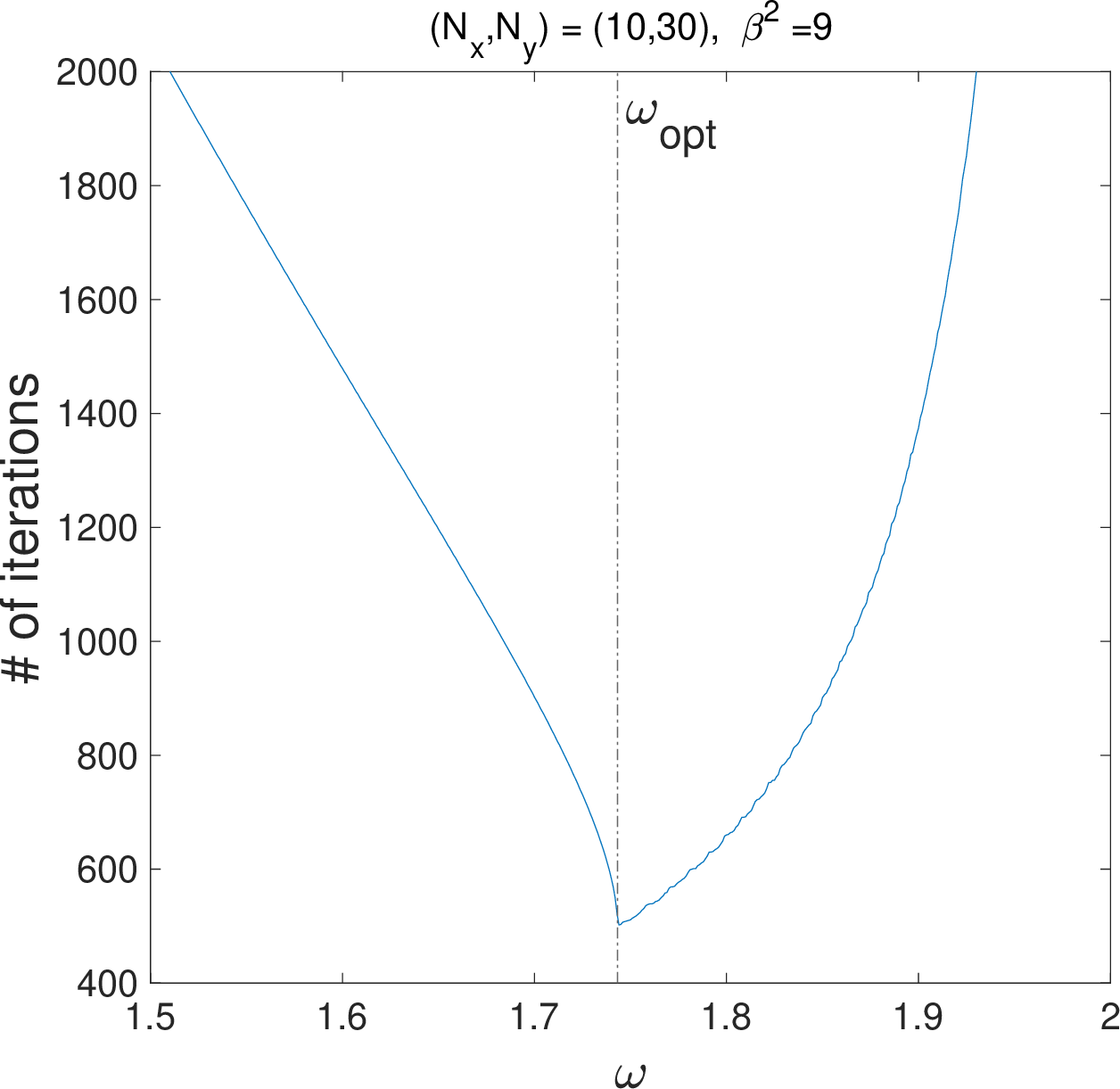}
	\includegraphics[width=0.5\textwidth]{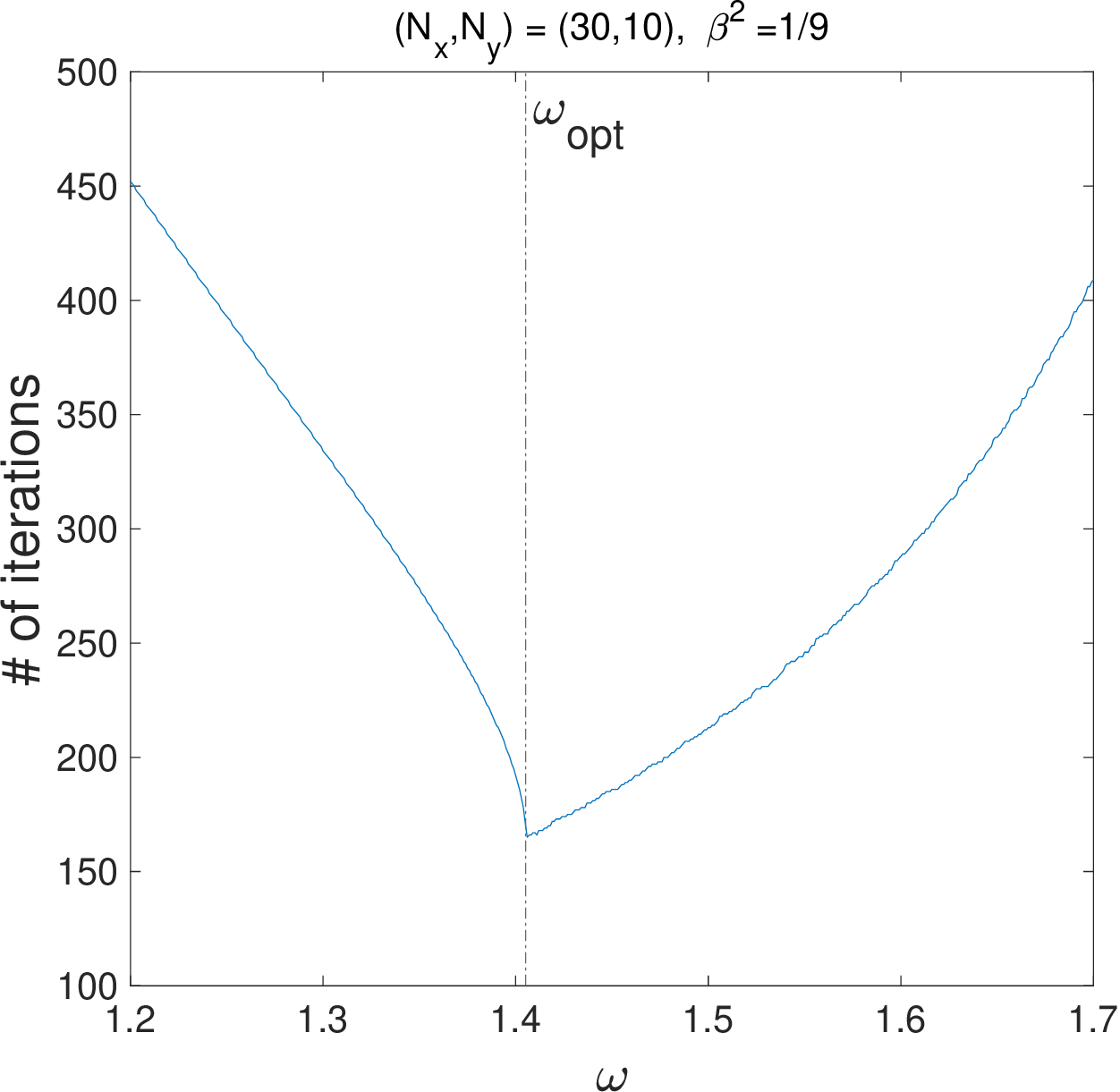}
	\caption{Number of iterations for the line SOR method of the HOC scheme \refp{e:lsor:hoc} for different values of $\omega$. The dash-dotted line indicates the optimal value \refp{e:opt:lsor:hoc}.}
	\label{fig:plot:iter:lsor:hoc}       
\end{figure}

Figure \ref{fig:plot:iter:lsor:hoc} shows the number of iterations required for the convergence. The results are very similar to those of the second-order scheme in Fig. \ref{fig:plot:iter:lsor}. Again, we observe the optimal relaxation parameter gives the least number of iterations for both the grids.  

The point SOR method for the HOC scheme is analyzed in \cite{adams1988sor} for $\beta=1$. Here, we follow the same approach. Separating the variables for \refp{e:compact:sor} leads to
\begin{align}
	\frac{\alpha^2 + \omega - 1 }{\omega} &= \frac{5-\beta^2}{10(1+\beta^2)}\frac{\alpha^2 X_{i-1} + X_{i+1}}{X_i} + \frac{5\beta^2-1}{10(1+\beta^2)}\frac{\alpha^2 Y_{j-1} + Y_{j+1}}{Y_j}   \nn\\
	&+ \frac{1}{20}
	\left(\frac{\alpha^2 Y_{j-1} + Y_{j+1}}{Y_j}\right)
	\left(\frac{X_{i-1} + X_{i+1}}{X_i}\right)
	\label{e:sep:compact}
\end{align}
Assuming the same function as \refp{e:YX} for $Y(y)$, we have
\begin{eqnarray}
	\label{e:eq1}
	\frac{\alpha^2 + \omega - 1 }{\omega} =
	\frac{5\beta^2-1}{5(1+\beta^2)}\alpha\cos(\theta_y) + \frac{X_{i-1}}{X_i}\Phi_1 + \frac{X_{i+1}}{X_i}\Phi_2
\end{eqnarray}
where
\begin{align*}
	\Phi_1 &= \frac{\alpha}{10}\cos(\theta_y) + \frac{5-\beta^2}{10(1+\beta^2)}\alpha^2 \\
	\Phi_2 &= \frac{\alpha}{10}\cos(\theta_y) + \frac{5-\beta^2}{10(1+\beta^2)}
\end{align*}

Now, assuming
\begin{equation}
	\label{e:hoc:X}
	X = \left(\frac{\Phi_1}{\Phi_2}\right)^{\frac{x}{2\Delta x}} \sin(k_x x)
\end{equation}
results in
\begin{eqnarray*}
	\frac{\Phi_1 X_{i-1} + \Phi_2 X_{i+1}}{X_i} = 2\sqrt{\Phi_1\Phi_2}\cos(k_x \Delta x)
\end{eqnarray*}
therefore, substituting this into \refp{e:eq1}, we have
\begin{eqnarray*}
	\frac{\alpha^2 + \omega - 1 }{\omega} =
	\frac{5\beta^2-1}{5(1+\beta^2)}\alpha\cos(\theta_y) + 2\sqrt{\Phi_1\Phi_2}\cos(\theta_x)
\end{eqnarray*}
and then by rearranging the terms, the following quartic equation is obtained
\begin{align}
	\alpha^4 &-\frac{2}{5}q\omega\left(c_1p^2\omega+5c_2\right)\alpha^3 
	+\left[ 2(\omega-1) + \left(c_2^2q^2-4c_1^2p^2-\frac{p^2q^2}{25} \right)\omega^2 \right]\alpha^2 \nn\\
	\label{e:quartic}
	&-\frac{2}{5}q\omega\left(c_1p^2\omega+5c_2(\omega-1)\right)\alpha + (\omega-1)^2 = 0 \\
	& c_1 = \frac{5-\beta^2}{10(1+\beta^2)} \; , \quad c_2 = \frac{4}{5}-2c_1 \; ,\qquad
	p=\cos(\theta_x)  \; , \quad q= \cos(\theta_y) \nn
\end{align}

Note that since $\beta^2 > 0$, we always have
\begin{eqnarray*}
	-\frac{1}{10}<c_1<\frac{1}{2}
\end{eqnarray*}

\begin{figure}[h]
	\setkeys{Gin}{draft=\showeps}
	\includegraphics[width=0.5\textwidth]{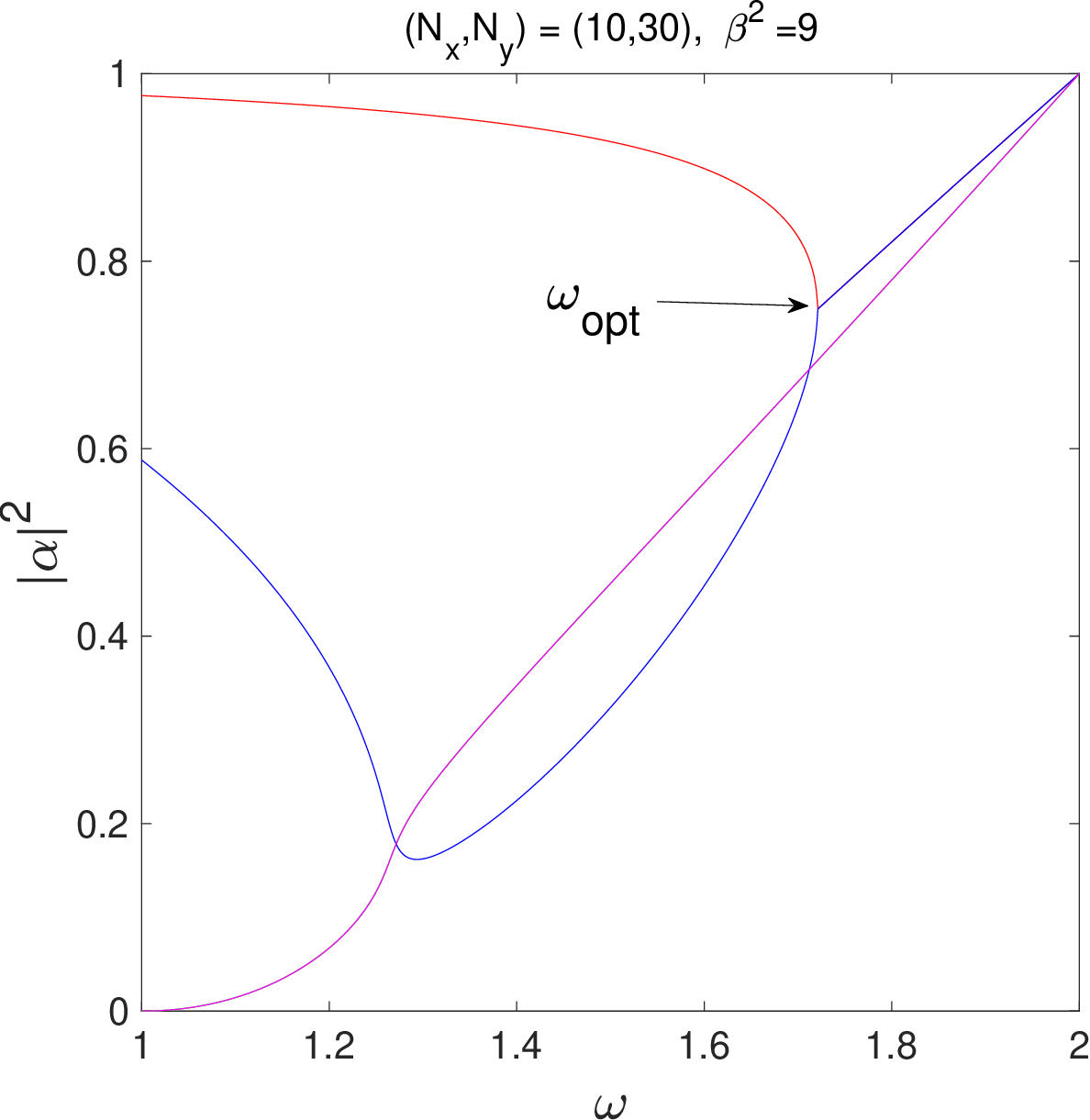}
	\includegraphics[width=0.5\textwidth]{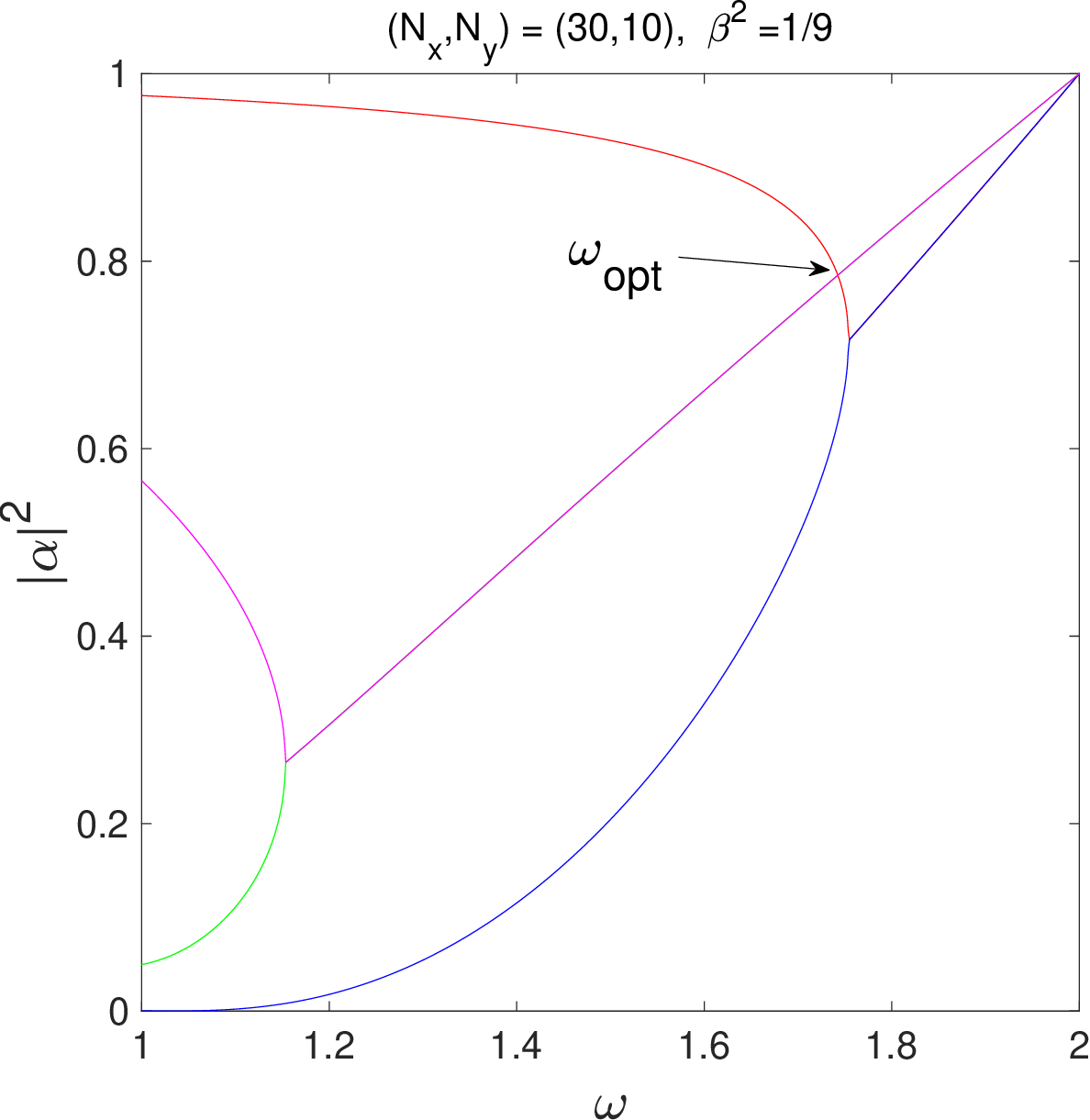}
	\caption{Distribution of $|\lambda| = |\alpha|^2$ obtained by solving \refp{e:quartic} for $k_x=k_y=\pi$ and $\omega\in [1,2]$.}
	\label{fig:plot}       
\end{figure}

For a fixed $\omega$, $k_x$ and $k_y$ the roots of \refp{e:quartic} can be obtained exactly. The roots may be real or complex. Figure \ref{fig:plot} shows the modulus of the eigenvalue ($|\lambda| = |\alpha|^2$) for $\omega\in[1,2]$ and $k_x=k_y=\pi$. The optimal $\omega$ occurs when these complex roots intersect the largest real root.  

Now, we try to find the optimal $\omega$ for small $h$ using the perturbation analysis. First, we perform a first-order analysis and then to improve the accuracy a second-order analysis is performed. For the first-order analysis we assume
\begin{eqnarray}
	\omega = 2 - k_1 h + O(h^2), \qquad h=\Delta x =\beta\Delta y
	\label{e:pert:w}
\end{eqnarray}

For now, we assume the largest root corresponds to $p=q=1$. To determine the limiting behavior we set $h=0$
\begin{align}
	\alpha^4 &-\frac{4}{5}q\left(2c_1p^2+5c_2\right)\alpha^3 
	+\left[ 2 + 4\left(c_2^2q^2-4c_1^2p^2-\frac{p^2q^2}{25} \right)\right]\alpha^2 \nn\\
	&-\frac{4}{5}q\left(2c_1p^2+5c_2\right)\alpha + 1 = 0 
	\label{e:quartic2}
\end{align}

In this limit all the roots lie on the unit circle (in the complex plane) since they must have modulus no greater than one (due to stability) and their product equals to $(\omega-1)^2=1$.
Choosing $\alpha=e^{i\phi}$ and dividing by $\alpha^2=e^{2i\phi}$, we have
\begin{eqnarray*}
	2\cos^2\phi -\frac{4}{5}q\left(2c_1p^2+5c_2\right)\cos\phi 
	+ 2\left(c_2^2q^2-4c_1^2p^2-\frac{p^2q^2}{25} \right) = 0 
	\label{e:cosine}
\end{eqnarray*}
where the roots are obtained as
\begin{align}
	\label{e:cosine1}
	\cos\phi_1 &= \frac{q}{5}\left(2c_1p^2+5c_2\right) - \frac{p}{5}\sqrt{\Delta'} \\
	\label{e:cosine2}
	\cos\phi_2 &= \frac{q}{5}\left(2c_1p^2+5c_2\right) + \frac{p}{5}\sqrt{\Delta'} \\
	\Delta'&= (1+20c_1c_2)q^2+4c_1^2p^2q^2+100c_1^2\nn
\end{align}
where $\Delta'$ can be shown that is always positive. For the first root, for all the values of $p,q\in [0,1]$ and $c_1\in [-0.1,0.5]$, we have $-1<\cos\phi_1<1$. Therefore, \refp{e:quartic2} always has at least one pair of non-real roots. For the second root, the right-hand side of \refp{e:cosine2} is a strictly increasing non-negative function of $p$ and $q$ for $p,q\in [0,1]$ and gives $\cos\phi_2=1$ only for $p=q=1$. Thus for $p=q=1$, there are two
real roots at $\alpha=1$, but for all other choices of $p$ and $q$, equation \refp{e:quartic2} has four non-real roots. We now wish to determine the behavior of the roots of \refp{e:quartic} as $h\rightarrow 0$. First, we consider the case of four non-real roots.

\subsubsection{Four non-real roots}
We assume
\begin{eqnarray}
	\label{e:betaroots}
	\alpha_{1,2}=e^{\pm i\phi_1}(1-\beta_1h)+O(h^2),\qquad \alpha_{3,4}=e^{\pm i\phi_2}(1-\beta_2h)+O(h^2)
\end{eqnarray}
since the product of the roots is $(1-\omega)^2$, equating the $O(h)$ terms yields
\begin{eqnarray}
	\label{e:rootMul}
	\beta_1+\beta_2=k_1
\end{eqnarray}
Also, the sum of the roots equals the coefficient of $\alpha^3$, i.e. 
\begin{eqnarray*}
	\frac{2}{5}q\omega (c_1p^2\omega-10c_1+4)
\end{eqnarray*}
and equating the $O(h)$ terms yields
\begin{eqnarray}
	\label{e:rootSum}
	-2\beta_1\cos\phi_1-2\beta_2\cos\phi_2 = -\frac{4}{5}k_1 q(2c_1p^2-5c_1+2) 
\end{eqnarray}
solving the equations \refp{e:rootMul} and \refp{e:rootSum}, gives
\begin{eqnarray}
	\label{e:beta12}
	\beta_1 = \left(\frac{1}{2}-\frac{c_1pq}{\sqrt{\Delta'}}\right)k_1 \; , \qquad
	\beta_2 = \left(\frac{1}{2}+\frac{c_1pq}{\sqrt{\Delta'}}\right)k_1 
\end{eqnarray}

The largest root of \refp{e:betaroots} corresponds to the smallest of $\beta_1$ and $\beta_2$, which depends on the sign of $c_1$. In any case we have 
\begin{eqnarray}
	\label{e:betamax}
	\beta_m = \min(\beta_1,\beta_2)=
	\left(\frac{1}{2}-\frac{|c_1|pq}{\sqrt{\Delta'}}\right)k_1 
\end{eqnarray}

The maximum of $\beta_m$ occurs for $p=q=1$, thus
\begin{eqnarray}
	\label{e:betamax2}
	\beta_m = \left(\frac{1}{2}-\frac{|c_1|}{1+8c_1}\right)k_1 
\end{eqnarray}

\subsubsection{Two real and two non-real roots}
For the real roots, we take terms out to the $O(h^2)$ term because now the $O(h)$ terms will cancel identically
\begin{eqnarray}
	\omega = 2 - k_1 h - k_2 h^2+O(h^3)\;, \qquad \alpha = 1-s_1h-s_2h^2+O(h^3)
	\label{e:pert:w2}
\end{eqnarray}
and also we replace the cosine functions with their approximations
\begin{align*}
	\cos\theta_x &= \cos(k_x \Delta x) = 1-\frac{1}{2}k_x^2 h^2 + O(h^4) \\
	\cos\theta_y &= \cos(k_y \Delta y) = 1-\frac{1}{2}\frac{k_y^2 h^2}{\beta^2} + O(h^4)
\end{align*}
then an expression for $s_1$ is obtained by equating the coefficients of $h^2$
\begin{eqnarray*}
	\frac{4}{5}(1+8c_1)s_1^2 - \frac{4}{5}(1+10c_1)k_1s_1 + 
	\frac{4}{25}(1+10c_1)^2(k_x^2+k_y^2)=0
\end{eqnarray*}
which gives
\begin{eqnarray*}
	s_1 = \frac{1+10c_1}{2(1+8c_1)}\left( k_1 \pm \sqrt{\Delta''}\right),\quad \Delta'' = k_1^2-\delta(k_x,k_y)^2
\end{eqnarray*}
where
\begin{eqnarray*}
	\delta(k_x,k_y) = \sqrt{\frac{4}{5}(1+8c_1)(k_x^2+k_y^2)}
\end{eqnarray*}

Equation \refp{e:quartic} has real roots if $s_1$ is real or $\Delta''$ is positive which leads to 
\begin{eqnarray*}
	k_1^2 \geq \delta(k_x,k_y)^2 = \frac{4}{5}(1+8c_1)(k_x^2+k_y^2)
\end{eqnarray*}
the largest real root occurs for the negative sign and by choosing $k_x=k_y=\pi$
\begin{eqnarray}
	\label{e:sm}
	s_m = \frac{1+10c_1}{2(1+8c_1)}\left( k_1 - \sqrt{k_1^2-\delta(\pi,\pi)^2}\right)
\end{eqnarray}
where $s_m$ is real for
\begin{eqnarray*}
	k_1 \geq \delta(\pi,\pi) = \pi\sqrt{\frac{8}{5}(1+8c_1)}
\end{eqnarray*}

We must still determine the other two roots which are complex. 
For these roots, we assume
\begin{eqnarray*}
	e^{\pm i\phi_3}(1-\beta_3h + O(h^2))
\end{eqnarray*}
since the product of the roots is $(1-\omega)^2$, equating the $O(h)$ terms yields
\begin{eqnarray}
	\beta_3 = \frac{1+6c_1}{2(1+8c_1)}k_1
\end{eqnarray}
which corresponds to $\beta_1$ in \refp{e:beta12}. This shows that for $c_1\geq 0$ ($\beta^2\leq 5$) it is the pair of complex roots with smaller modulus that splits into real roots. For $c_1<0$ this occurs for the pair of complex roots with larger modulus. This affects the optimal relaxation factor position as seen in Fig. \ref{fig:plot}. To find the optimal $\omega$, we require to minimize the following expression
\begin{eqnarray*}
	\max(1-\beta_m h , 1-s_mh) 
\end{eqnarray*}
since $\beta_m$ is an increasing function and $s_m$ is a decreasing function of $k_1$, the minimum occurs for $s_m=\beta_m$ which results in
\begin{eqnarray}
	\label{e:k1}
	k_{1,opt} = \delta(\pi,\pi)\times
	\left\{
	\begin{array}{lr}
		\frac{1+10c_1}{\sqrt{(1+14c_1)(1+6c_1)}}\quad
		& c_1\geq 0 \\&\\	
		1
		& c_1<0
	\end{array}
	\right.
\end{eqnarray}
therefore
\begin{eqnarray}
	\label{e:wopt:hoc1}
	\omega_{opt} = 2-k_{1,opt} h + O(h^2)
\end{eqnarray}

Figure \ref{fig:plot:sor:hoc} shows the position of this approximation (shown by $\omega_1$). The figure shows the first-order approximation is not accurate enough.  
We try to find a more accurate approximation by considering the second-order term. 

\subsubsection{Second-order perturbation analysis}
For the second-order perturbation analysis, we again consider both the cases of the roots considered in the first-order analysis. Now, we assume
\begin{eqnarray}
	\omega = 2 - k_1 h-k_2 h^2 + O(h^3)
	\label{e:pert:w_2}
\end{eqnarray} 
and also for the four complex roots, we assume  
\begin{eqnarray}
	\label{e:betaroots_2}
	e^{\pm i\phi_1}(1-\beta_1 h-\gamma_1 h^2)+O(h^3),\qquad e^{\pm i\phi_2}(1-\beta_2 h-\gamma_2 h^2)+O(h^3)
\end{eqnarray}
equating the coefficients of $h^2$ for the roots product and roots sum, we obtain
\begin{align*}
	&\gamma_1+\gamma_2=k_2+Bk_1^2 \;, \qquad B =  \frac{1}{4}-\frac{{c_{1}}^2\,p^2\,q^2}{\Delta'}\\
	&\gamma_1\cos\phi_1+\gamma_2\cos\phi_2 = c_2k_2q + \frac{4k_2-k_1^2}{5}c_1p^2q 
\end{align*}
where using \refp{e:cosine1} and \refp{e:cosine2} gives
\begin{align}
	\label{e:gamma1}	    
	\gamma_1&= \left(\frac{1}{2}-\frac{c_1pq}{\sqrt{\Delta'}}\right)k_2+\frac{1}{2}\left(B+\frac{(2c_{1}p^2+5c_{2})Bq+c_1p^2q}{p\sqrt{\Delta'}}\right)k_1^2 \\
	\label{e:gamma2}	    
	\gamma_2&= \left(\frac{1}{2}+\frac{c_1pq}{\sqrt{\Delta'}}\right)k_2+\frac{1}{2}\left(B-\frac{(2c_{1}p^2+5c_{2})Bq+c_1p^2q}{p\sqrt{\Delta'}}\right)k_1^2
\end{align}
the value which corresponds to $\beta_m$ \refp{e:betamax2} can be expressed as
\begin{eqnarray}
	\label{e:gammam}
	\gamma_m =  A k_2+D k_1^2
\end{eqnarray}
where
\begin{eqnarray}
	A = \frac{1}{2}-\frac{|c_1|pq}{\sqrt{\Delta'}}\;,\qquad
	D = \frac{1}{2}\left(B+\mathrm{sign}(c_1)\frac{(2c_{1}p^2+5c_{2})Bq+c_1p^2q}{p\sqrt{\Delta'}}\right) \nonumber	
\end{eqnarray}

For the real roots we require the third-order terms
\begin{align}
	\omega &= 2 - k_1 h - k_2 h^2 - k_3h^3+O(h^4)\nonumber\\
	\alpha &= 1-s_1h-s_2h^2-s_3h^3+O(h^4)
	\label{e:pert:w3}
\end{align}
where by inserting in \refp{e:quartic} and equating the coefficient of $h^3$ to zero, we obtain
\begin{eqnarray}
	s_2 
	= \frac{s_1}{2}\left(-s_1+\frac{k_1^2+2k_2}{\pm\sqrt{\Delta''}}\right)
	\label{e:s2}
\end{eqnarray}
the value of $s_2$ which corresponds to $s_m$ \refp{e:sm}, by setting $k_x = k_y = \pi$, is given by 
\begin{eqnarray}
	\label{e:s2m}
	s_{2,m} 
	= \frac{s_m}{2}\left(-s_m+\frac{k_1^2+2k_2}{-R_\pi}\right)
	\;,\qquad 	R_\pi = \sqrt{k_1^2-\delta(\pi,\pi)^2}
\end{eqnarray}
now equating $\gamma_m = s_{2,m}$ gives
\begin{eqnarray}
	\label{e:k2}
	k_{2,opt} = 
	-\frac{k_1^2}{2A}\times\frac{R_\pi A^2+A k_1+2R_\pi D}{R_\pi + k_1}
\end{eqnarray}
therefore
\begin{eqnarray}
	\label{e:wopt:hoc2}
	\omega_{opt} = 2-k_{1,opt} h - k_{2,opt} h^2 + O(h^3)
\end{eqnarray}

\begin{figure}[h]
	\setkeys{Gin}{draft=\showeps}
	\includegraphics[width=0.5\textwidth]{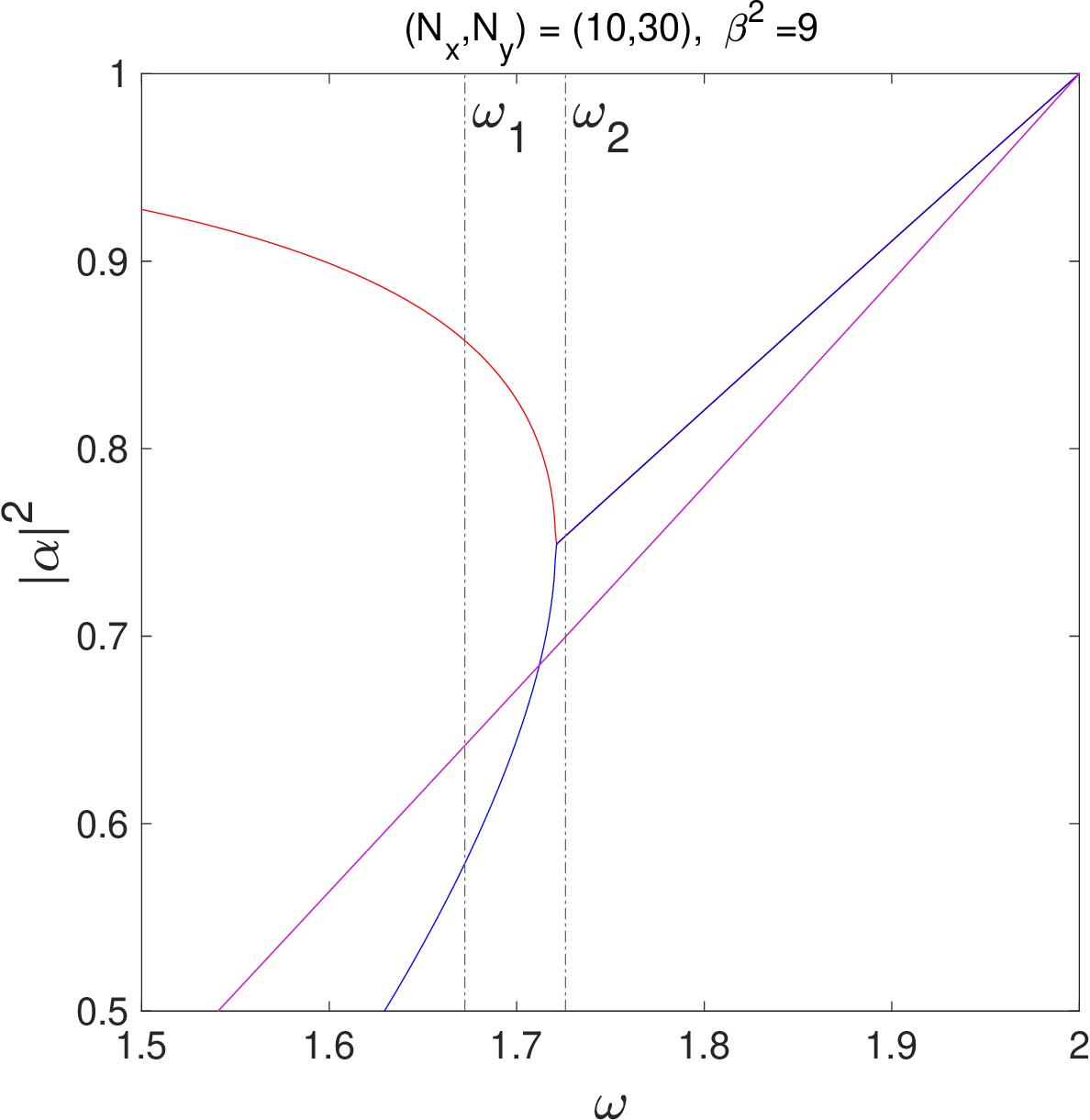}
	\includegraphics[width=0.5\textwidth]{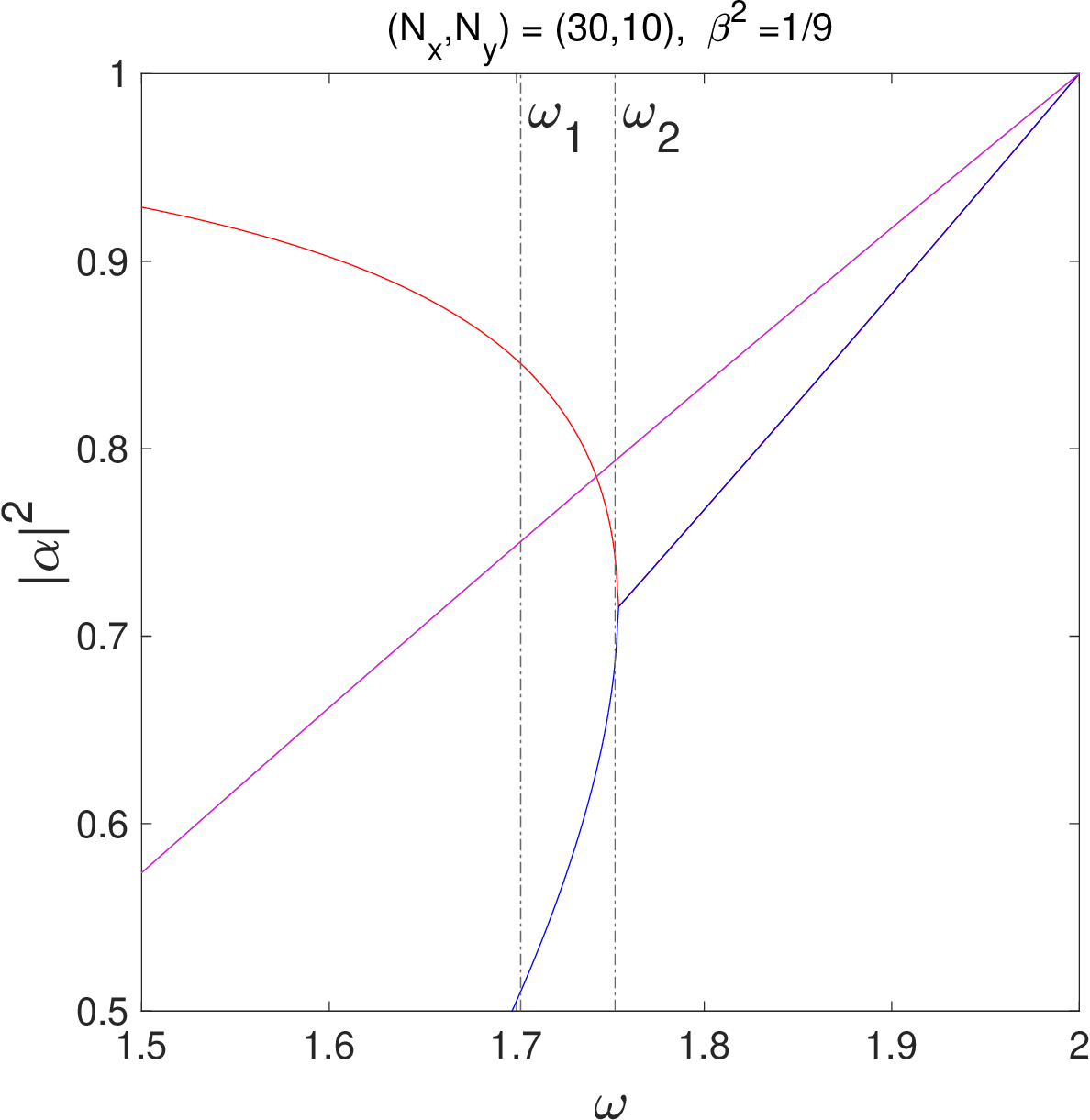}
	\caption{The optimal relaxation factor for the point SOR method of the HOC scheme using first-order \refp{e:wopt:hoc1} and second-order \refp{e:wopt:hoc2} approximations, shown by $\omega_1$ and $\omega_2$, respectively.}
	\label{fig:plot:sor:hoc} 
\end{figure}

Figure \ref{fig:plot:sor:hoc} shows the position of this approximation (shown by $\omega_2$). The figure shows the second-order approximation is much more accurate than the first-order one.

\begin{figure}[h]
	\setkeys{Gin}{draft=\showeps}
	\includegraphics[width=0.45\textwidth]{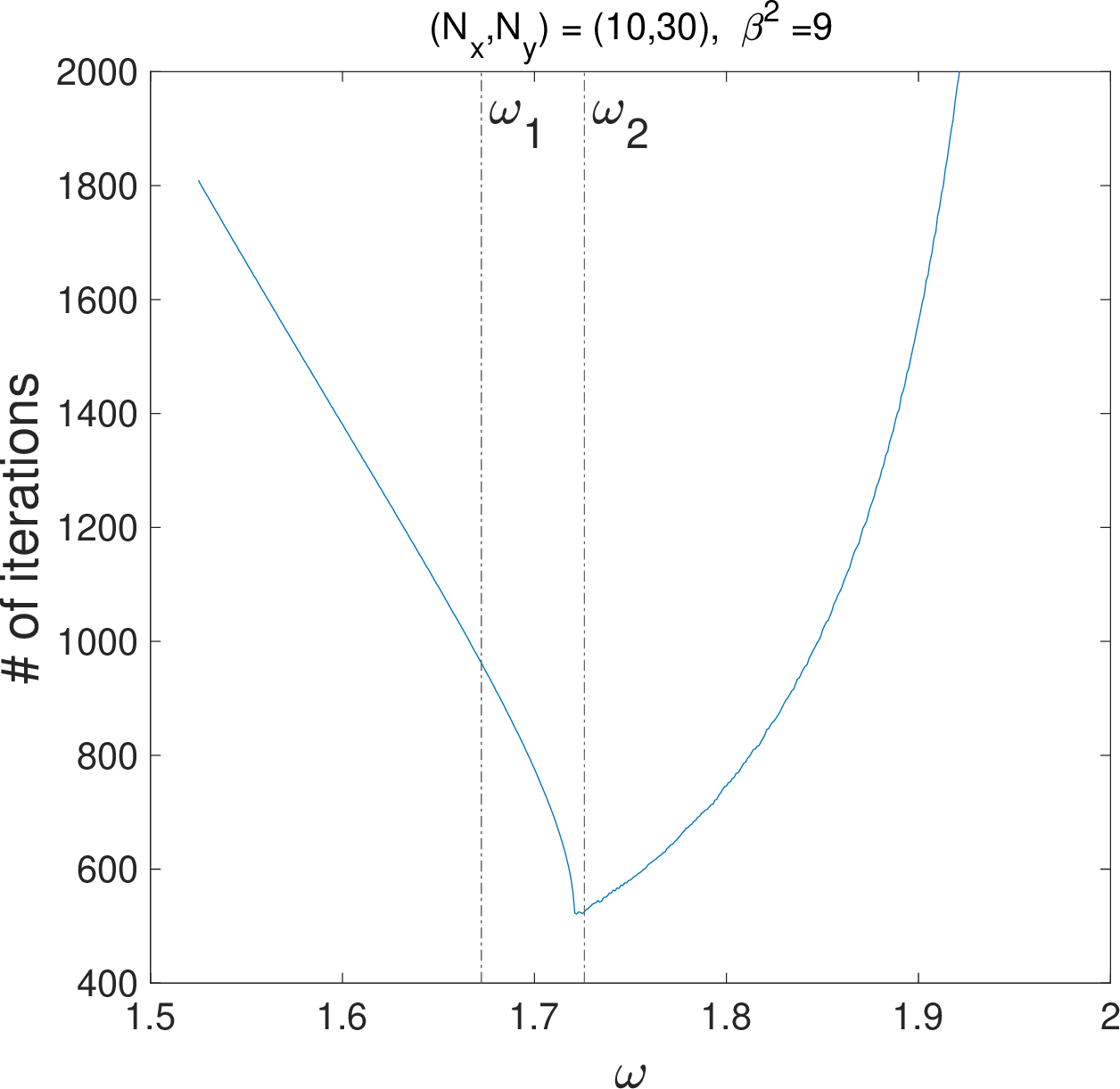}
	\includegraphics[width=0.45\textwidth]{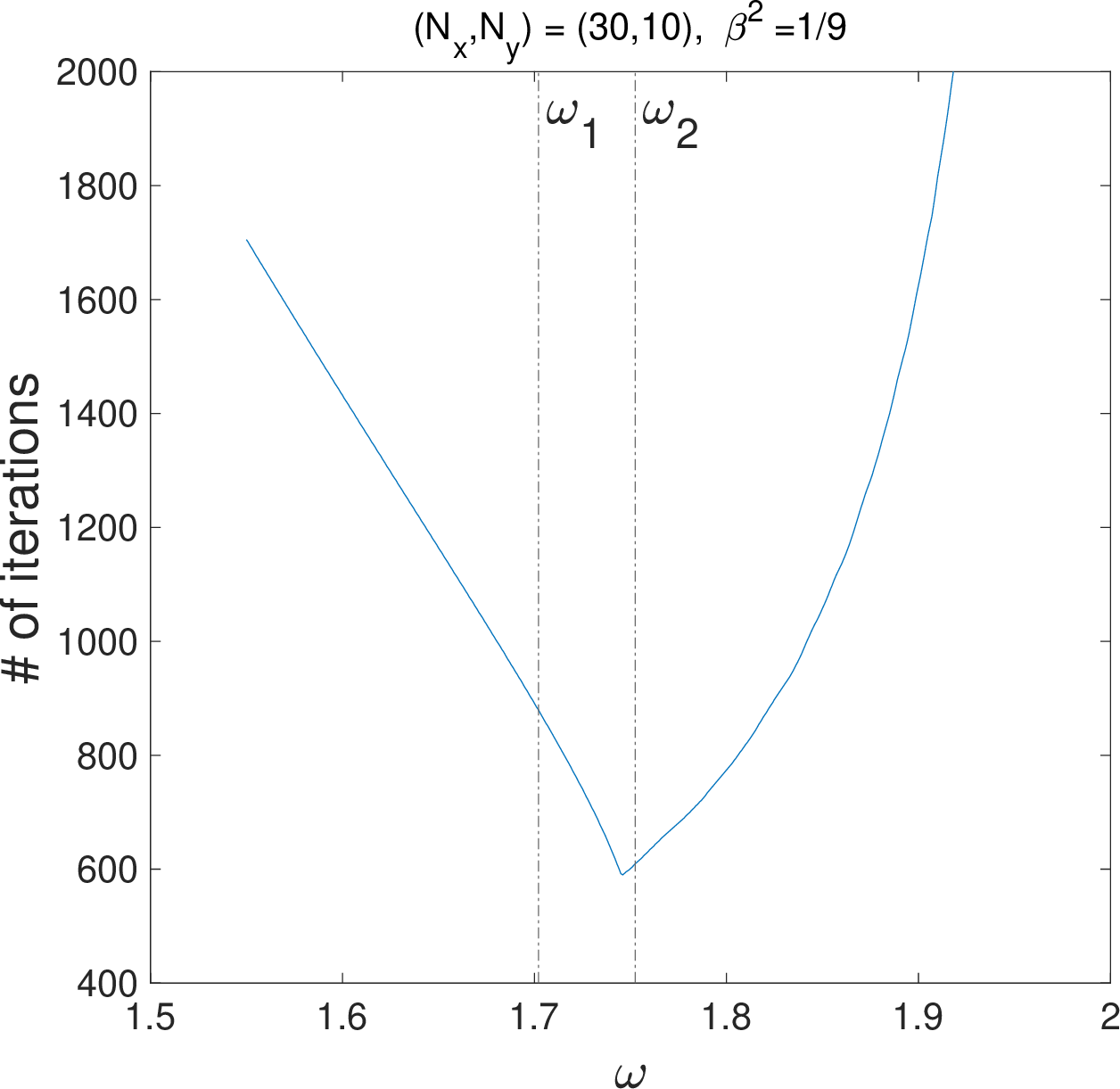}
	\caption{Number of iterations for the point SOR method of the HOC scheme \refp{e:compact:sor} for different values of $\omega$. The dash-dotted lines indicate the first-order \refp{e:wopt:hoc1} and second-order \refp{e:wopt:hoc2} approximations, shown by $\omega_1$ and $\omega_2$, respectively.}
	\label{fig:plot:iter:sor:hoc}
\end{figure}

Figure \ref{fig:plot:iter:sor:hoc} shows the number of iterations required for the convergence for different values of $\omega$. The figure shows the  number of iterations for the first-order approximation \refp{e:wopt:hoc1} is nearly $25\%$ more than that of the actual optimal relaxation parameter. Also, we observe the number of iterations for the second-order approximation \refp{e:wopt:hoc2} is almost the same as those of the actual optimal relaxation parameter. 

\section{Neumann boundary conditions}
\label{sec:4}
The analysis conducted in the previous sections assumes Dirichlet boundary conditions for all the boundaries. Here, we obtained the optimal relaxation factor when one or more boundaries are of Neumann type. First, let the normal derivative is known on the right edge:
\begin{eqnarray}
	\label{e:bcR}
	\left.\frac{\partial u}{\partial x}\right|_{x=1}=0
\end{eqnarray}
and all the other boundaries are of Dirichlet type. In this case, we extend the domain by mirroring against the right edge. Therefore, the new domain is $\Omega = [0,2]\times [0,1]$. Now, if all the boundaries of this domain are of Dirichlet type ($u=0$), the solution will be symmetric around $x=1$:
\begin{eqnarray}
	u(x,y)=u(2-x,y)
	\label{e:symmetry}
\end{eqnarray}

This symmetry causes the boundary condition \refp{e:bcR} to be automatically satisfied which means the optimal relaxation parameter can be obtained by \refp{e:opt:2nd}. However, due to extension of the domain, the number of grid points is doubled along $x$-direction ($2N_x+1$ points) and the allowable values for $k_x$ (equation \refp{e:kx}) are changed as
\begin{eqnarray*}
	k_x = \frac{1}{2} n_x \pi\;,\quad 1 \leq  n_x\leq {2N_x-1}
\end{eqnarray*}
which means the largest possible value of $r$ is obtained by $k_x = \frac{\pi}{2}$ instead of $k_x = \pi$. Therefore, we have
\begin{eqnarray}
	\omega_{opt} = \frac{2}{1 + \sqrt{1-r_{opt}^2}}, \qquad
	r_{opt} = \frac{\cos(\frac{1}{2}\pi\Delta x) +\beta^2\cos(\pi\Delta y)}{1+\beta^2}	
	\label{e:opt:2nd:neumannX}
\end{eqnarray}

If the symmetry is applied by using the ghost-node method in Equations \refp{e:central:sor}-\refp{e:lsor:hoc} for the right edge points of the original domain ($i=N_x$)
\begin{eqnarray*}
	u_{N_x+1,j}=u_{N_x-1,j}
\end{eqnarray*}
then \refp{e:opt:2nd:neumannX} gives the optimal relaxation factor for this case. 

Now, let both the left and right boundaries are of Neumann type:
\begin{eqnarray}
	\label{e:bcRL}
	\left.\frac{\partial u}{\partial x}\right|_{x=0}=0\;,\qquad \left.\frac{\partial u}{\partial x}\right|_{x=1}=0
\end{eqnarray}
For this case, the assumption for $X$ in \refp{e:YX} does not satisfy the boundary conditions and the sine function should be replaced by cosine as
\begin{eqnarray}
	X = \alpha^{\frac{x}{\Delta x}} \cos(k_x x)
	\label{e:opt:2nd:XX}
\end{eqnarray}
Fortunately, this assumption again reduces to the same formulations which resulted in \refp{e:sor:r}. However, zero is added to the allowable values of $k_x$:   
\begin{eqnarray}
	k_x = n_x \pi\;,\quad 0 \leq  n_x\leq {N_x-1}
	\label{e:allowedcosx} 
\end{eqnarray}
This means the largest possible value of $r$ corresponds to $\theta_x = 0$ and $\theta_y = \pi \Delta y$. Hence,  
\begin{eqnarray}
	\omega_{opt} = \frac{2}{1 + \sqrt{1-r_{opt}^2}}, \qquad
	r_{opt} = \frac{1 +\beta^2\cos(\pi\Delta y)}{1+\beta^2}	
	\label{e:opt:2nd:neumannXX}
\end{eqnarray}

\begin{figure}[h]
	\setkeys{Gin}{draft=\showeps}
	\includegraphics[width=0.45\textwidth]{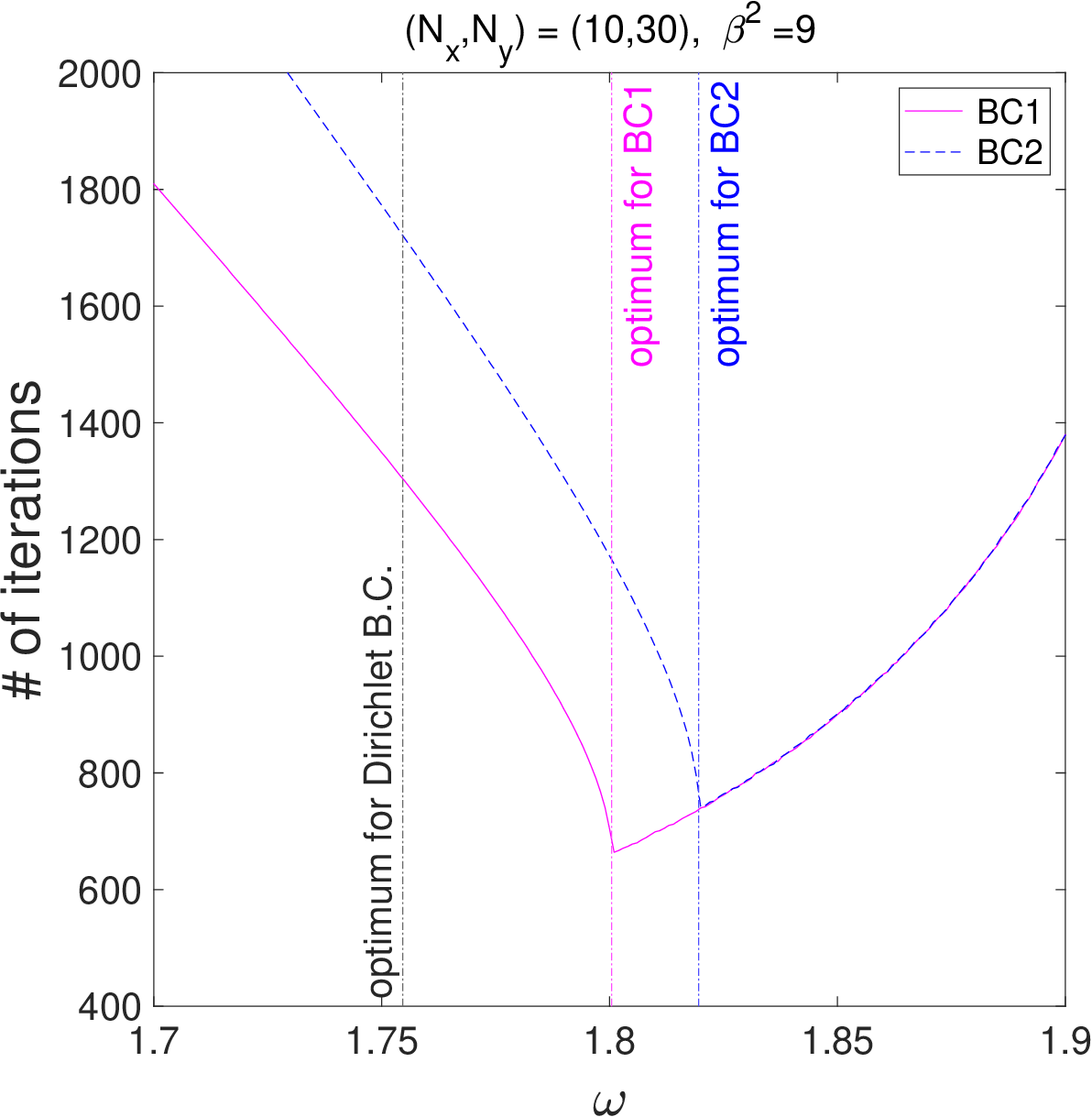}
	\includegraphics[width=0.45\textwidth]{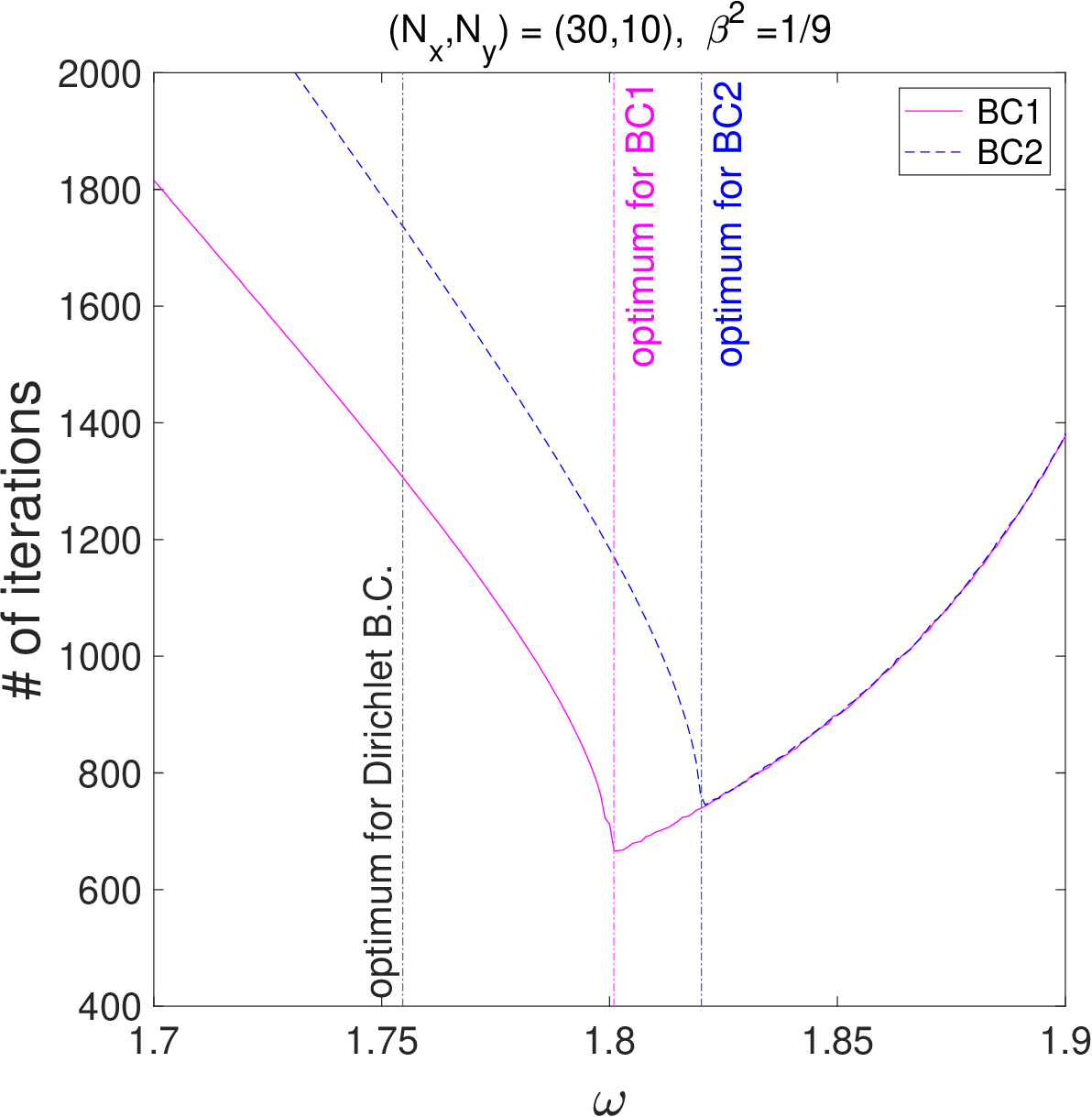}
	\caption{Number of iterations for the point SOR method of the second-order scheme \refp{e:central:sor} for different values of $\omega$. BC1: Neumann B.C. for the right edge. BC2: Neumann B.C. for both the left and right edges. The dash-dotted lines indicate the optimal values for each case (\refp{e:opt:2nd:neumannX} and \refp{e:opt:2nd:neumannXX}).}
	\label{fig:plot:iter:sor:bc}
\end{figure}

Figure \ref{fig:plot:iter:sor:bc} shows the number of iterations for the above mentioned boundary  conditions. The results verify the optimal values obtained by \refp{e:opt:2nd:neumannX} and \refp{e:opt:2nd:neumannXX}. Also, we observe the number of iterations using the optimal value for Dirichlet boundary conditions \refp{e:opt:2nd} is nearly twice more than that of the actual optimal relaxation parameter. 

\section{Robin boundary conditions}
\label{sec:5}
In this section, we consider the third type of boundary conditions or Robin boundary conditions which is specified by a linear combination of the function and its normal derivative on the boundary. For the $x$-direction the Robin boundary conditions are as follows 
\begin{eqnarray}
	\label{e:bcRobin}
	\left(au+b\frac{\partial u}{\partial x}\right)_{x=0}=0\;, \qquad
	\left(cu+d\frac{\partial u}{\partial x}\right)_{x=1}=0
\end{eqnarray}
where $a$, $b$, $c$ and $d$ are constant numbers and at least three of them are non-zero (otherwise \refp{e:bcRobin} reduces to Dirichlet or Neumann type). We need to choose a function for $X$ in \refp{e:YX} which simultaneously satisfies the discretized form of \refp{e:bcRobin}: 
\begin{eqnarray}
	\label{e:bcRobinDiscr}
	\left(aX_i+b\frac{X_{i+1}-X_{i-1}}{2\Delta x}\right)_{i=0}=0\;, \qquad
	\left(cX_i+d\frac{X_{i+1}-X_{i-1}}{2\Delta x}\right)_{i=N_x}=0
\end{eqnarray}
Therefore, we assume
\begin{eqnarray}
	X(x) = \alpha^{\frac{x}{\Delta x}} \sin(k_x x + \psi)\;,\qquad 
	\label{e:robin:X}
\end{eqnarray}
where again reduces to \refp{e:sor:r}. Now, this function should be replaced in \refp{e:bcRobinDiscr}. However, to reduce the analysis complexity, we ignore $\alpha^{\frac{x}{\Delta x}}$ variations, which is reasonable when $|\alpha|\simeq 1$. Using this approximation, we obtain
\begin{eqnarray}
	\label{e:psiL}
	&a\sin\psi+b\, \skx\cos\psi = 0\\
	\label{e:psiR}
	&c\sin(k_x+\psi)+d\, \skx\cos(k_x+\psi) = 0
\end{eqnarray}
where
\begin{eqnarray*}
	\skx = \frac{\sin(k_x \Delta x)}{\Delta x} = N_x \sin(\frac{k_x}{N_x})
\end{eqnarray*}

Expanding \refp{e:psiR} and then collecting $\sin\psi$ and $\cos\psi$ terms, give
\begin{eqnarray}
	\label{e:bcRobinL}
	&a\sin\psi+b\; \skx\cos\psi=0 \\
	\label{e:bcRobinR}
	&(c\cos k_x - d \,  \skx\sin k_x) \sin\psi + (c\sin k_x+d\, \skx\cos k_x) \cos\psi=0
\end{eqnarray}

Now, we investigate the allowable values for $k_x$. If $k_x=0$ is allowable, then we must have $\sin\psi=0$ or $a=c=0$. The former gives a zero function and the latter is equivalent to the case where both the left and right boundaries are of Neumann type which was presented previously. Therefore, we assume $k_x \neq 0$ from now on.

Equations \refp{e:bcRobinL} and \refp{e:bcRobinR} have a solution if
\begin{eqnarray*}
	\frac{a}{c\cos k_x - d\, \skx^2 \sin k_x}=\frac{b\, \skx^2}{c\sin k_x+d\, \skx^2\cos k_x}
\end{eqnarray*}
or, after rearranging
\begin{eqnarray}
	\label{e:bcRobin:kx}
	(ac+\skx^2bd )\sin k_x + (ad-bc)\, \skx\cos k_x=0
\end{eqnarray}

The smallest positive $k_x$ which satisfies \refp{e:bcRobin:kx} determines the optimal relaxation parameter. The left hand side of the above equation is an odd function of $k_x$ and also periodic with period $2\pi N_x$. Therefore, excluding $k_x=0$, searching in $(0,\pi N_x)$ is enough for finding all the roots. 

However, it is not clear how many roots are exactly in the interval. In fact, for some values of $a$, $b$, $c$ and $d$ the number of roots is less than the number of required eigenvalues. Therefore, we assume a hyperbolic function for $X(x)$:
\begin{eqnarray}
	X(x) = \alpha^{\frac{x}{\Delta x}} \sinh(k_x x + \psi)
	\label{e:robin:hyper:X}
\end{eqnarray}
which causes all the $\cos(\theta_x)$ terms be replaced by $\cosh(\theta_x)$ in the calculation of the eigenvalues. For instance, equation \refp{e:sor:r} changes to
\begin{eqnarray}
	r = \frac{\cosh(\theta_x) +\beta^2\cos(\theta_y)}{1+\beta^2}	
	\label{e:sor:hyper:r}
\end{eqnarray}
which means we may have $r>1$ and therefore some eigenvalues are greater than one in modulus. Using the similar calculations which resulted in \refp{e:bcRobin:kx}, we obtain the following equation
\begin{eqnarray}
	\label{e:bcRobin:hyper:kx}
	(ac-\shkx^2bd )\sinh k_x + (ad-bc)\, \shkx\cosh k_x=0
\end{eqnarray}
where
\begin{eqnarray*}
	\shkx = \frac{\sinh(k_x \Delta x)}{\Delta x} = N_x \sinh(\frac{k_x}{N_x})
\end{eqnarray*}

The left hand side of \refp{e:bcRobin:hyper:kx} is an odd function of $k_x$ and searching in $(0,+\infty)$ is enough for finding all the roots. The largest (non-zero) root of this equation, if any exists, determines the convergence behavior of the iteration method. Otherwise, the smallest root of \refp{e:bcRobin:kx} determines the optimal relaxation parameter.

In general, equations \refp{e:bcRobin:hyper:kx} and \refp{e:bcRobin:kx} must be solved numerically. However, some comments are given in the following. 

\subsection{$ad-bc=0$}
If $ad-bc=0$, which means none of the coefficients is zero (since we assumed three of the coefficients are non-zero), 
then equation \refp{e:bcRobin:hyper:kx} has a non-zero root at $k_x = \sinh^{-1}(\left|\frac{a}{b}\right|\Delta x)/\Delta x$.

\subsection{$ad-bc\neq0$}
Table \ref{t:mn} shows the number of positive roots for equation \refp{e:bcRobin:hyper:kx} based on the following parameters
\begin{eqnarray*}
	m = \frac{ac}{ad-bc}\;, \qquad n=\frac{bd}{ad-bc}
\end{eqnarray*}
where the details are given in \ref{s:app}. Note that, the definition of $m$ and $n$ limits their allowable values to $1+4mn\geq 0$.
\begin{table}[h]
	\begin{center}
		\caption{Number of positive roots for equation \refp{e:bcRobin:hyper:kx}.}
		\label{t:mn}
		\begin{tabular}{|c|c|c|c|}\hline
			& $m+1<0$ & $m+1=0$ & $m+1>0$  \\\hline
			$n>0$ &at most 2&at most 1&1\\\hline
			$n\leq 0$ &1&0&0\\\hline
		\end{tabular}
	\end{center}
\end{table}

For $m+1=0$, equation \refp{e:bcRobin:hyper:kx} has a double root at $k_x = 0$ (this is also the case for equation \refp{e:bcRobin:kx}). Therefore, for this case $k_x=0$ is the desired root if no positive root exists. 

\begin{figure}[h]
	\setkeys{Gin}{draft=\showeps}
	\includegraphics[width=0.45\textwidth]{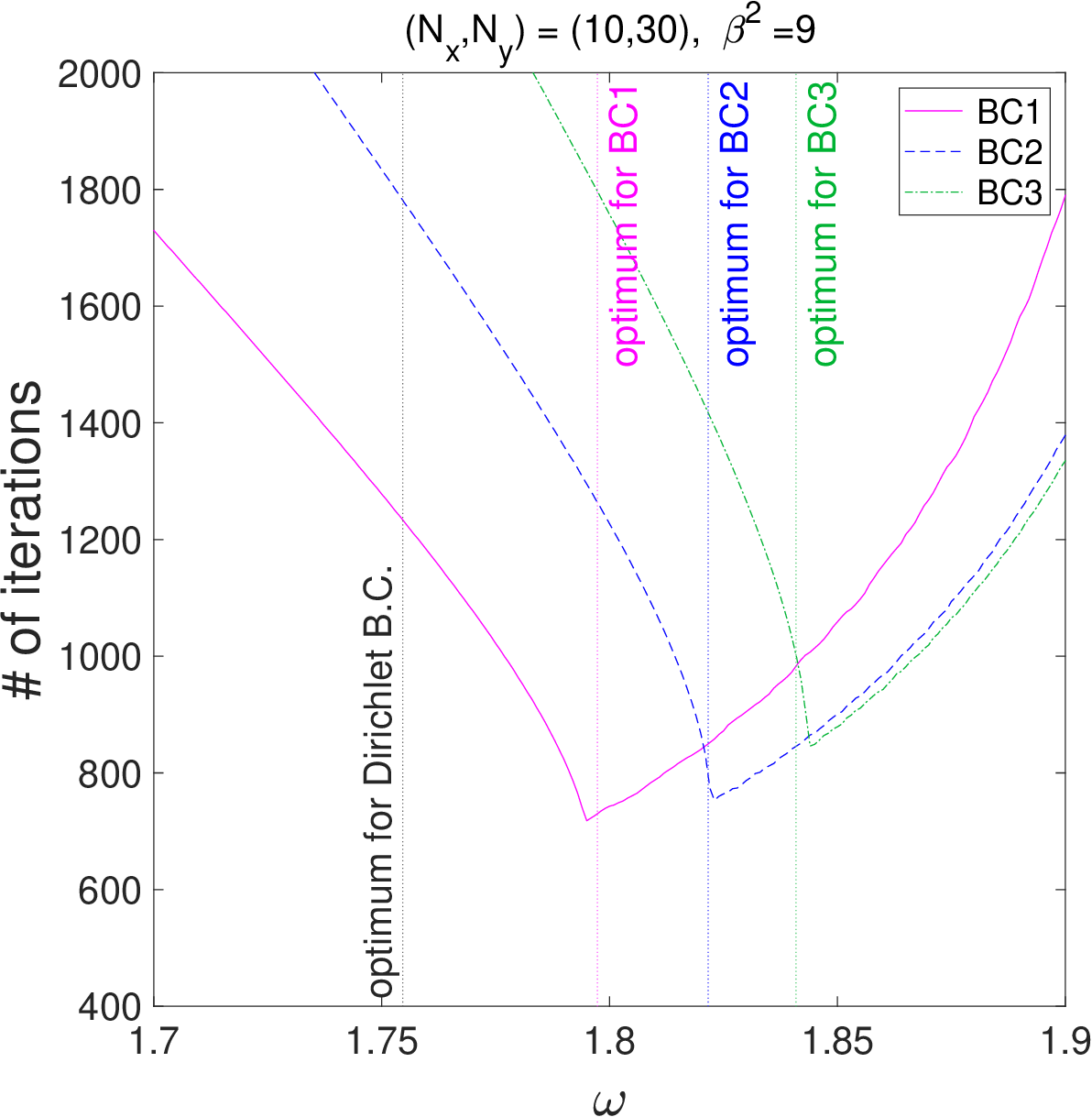}
	\includegraphics[width=0.45\textwidth]{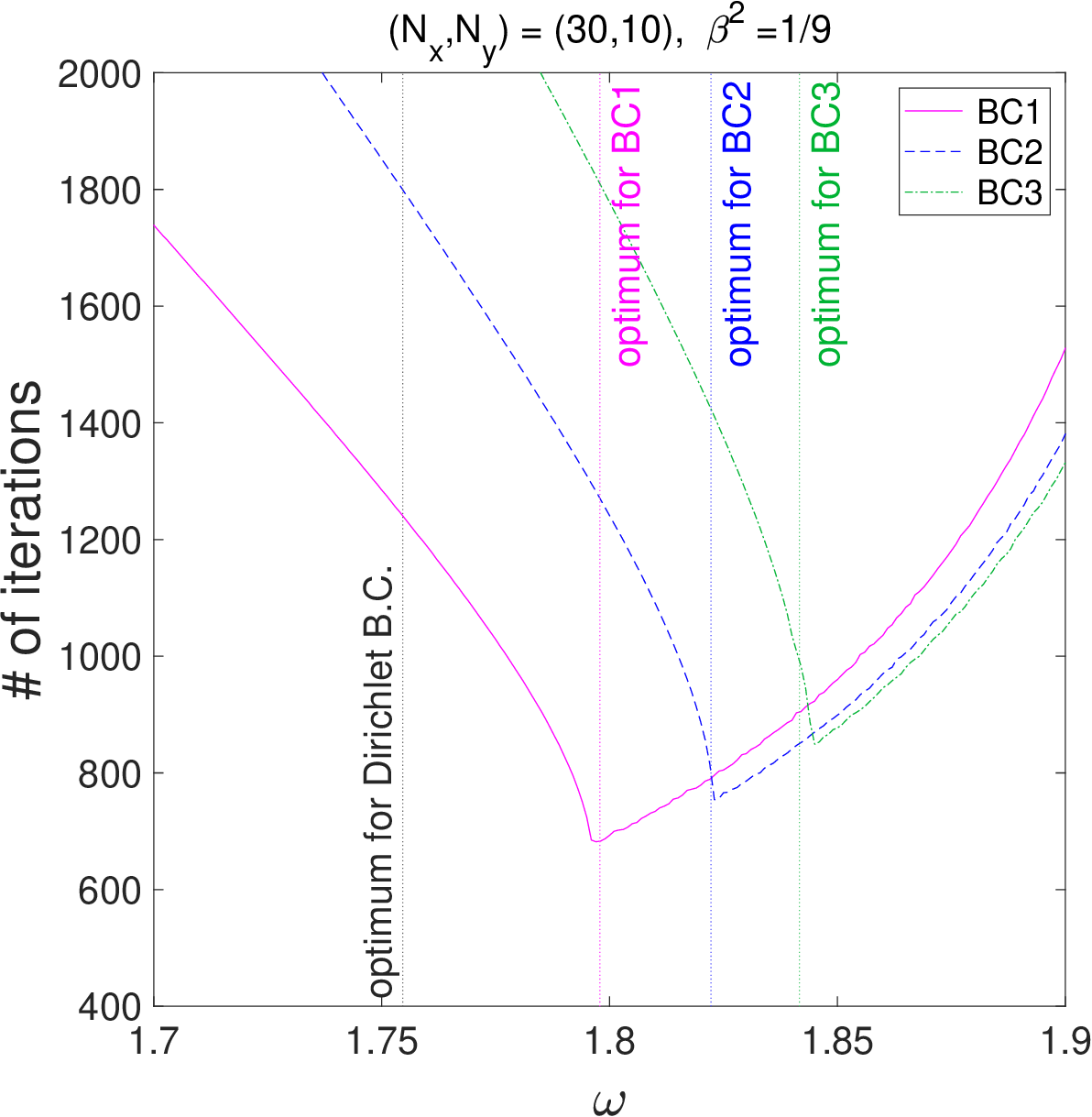}
	\caption{Number of iterations for the point SOR method of the second-order scheme \refp{e:central:sor} for different cases of the Robin boundary conditions (Table \ref{t:BCs}).}
	\label{fig:plot:iter:sor:robin}
\end{figure}

Figure \ref{fig:plot:iter:sor:robin} shows the number of iterations for several cases of the Robin boundary conditions where the coefficients are given in Table \ref{t:BCs}. We observe a slight difference between the analytical (vertical lines) and numerically found optimal parameters. This is expected, because we made an approximation by ignoring the variation of $\alpha^{\frac{x}{\Delta x}}$ in \refp{e:robin:X} when replaced in \refp{e:bcRobinDiscr}. However, the numerical results show the difference is very small.

\begin{table}[h]
	\begin{center}
		\caption{The coefficients of the Robin boundary conditions \refp{e:bcRobin} in Figures \ref{fig:plot:iter:sor:robin} and \ref{fig:plot:iter:hoc:robin}.}
		\label{t:BCs}
		\begin{tabular}{|c|rrrr|c|c|c|c|c|}
			\hline
			&$a$, &   $b$, & ~~~$c$, & ~~~$d$ & $ad-bc$ & $m$  & $n$  & equation & $k_x$ \\ \hline
			BC1&  1, & -0.25, &      1, &      1 &  1.25   & 0.8  & -0.2 & \refp{e:bcRobin:kx} & 1.70073 \\ \hline
			BC2&  1, &     2, &      1, &      2 &    0    & ---  & ---  & \refp{e:bcRobin:hyper:kx} & 0.49998 \\ \hline
			BC3&  1, &     1, &      1, &     -1 &   -2    & -0.5 & 0.5  & \refp{e:bcRobin:hyper:kx} & 1.54300 \\ \hline
		\end{tabular}
	\end{center}
\end{table}

The same analysis can be done for the $y$-direction boundary conditions (top and bottom edges). In summary, the following equation along with Table \ref{t:2ndbcsummary} can be used to obtain the optimal relaxation parameter for different cases of the boundary conditions. 
\begin{eqnarray}
	\omega_{opt} = \frac{2}{1 + \sqrt{1-r(k_x,k_y)^2}}, \qquad r(k_x,k_y) = \frac{\cos(k_x\Delta x) +\beta^2\cos(k_y\Delta y)}{1+\beta^2}
	\label{e:rfunc}
\end{eqnarray}

Note that if all the boundary conditions are of Neumann type, we obtain $r=1$ and $\omega_{opt}=2$ which is not convergent. This is not surprising since in this case the solution of the Poisson equation is not unique.
\begin{table}[h]
	\caption{The suitable values for $k_x$ and $k_y$ for different types of boundary conditions.} 
\label{t:2ndbcsummary}
\centering
\begin{tabular}{|c|c|}\hline
	$x$- or $y$-direction boundary types  &	$k_x$ or $k_y$ \\\hline
	Both Dirichlet  & $\pi$ \\\hline
	One Neumann & $\frac{1}{2}\pi$ \\\hline
	Both Neumann & $0$ \\\hline
	One Robin & The largest positive root of \refp{e:bcRobin:hyper:kx} \\
	& otherwise, the smallest root of \refp{e:bcRobin:kx} \\\hline
\end{tabular}
\end{table}


The optimal value for the other schemes can be computed using the same approach. 
For the line SOR schemes, \refp{e:opt:lsor} and \refp{e:opt:lsor:hoc} must be rewritten, respectively, as
\begin{eqnarray}
r(k_x,k_y) = \frac{\beta^2\cos(k_y\Delta y)}{1+\beta^2-\cos(k_x\Delta x)}
\label{e:opt:lsor:bc}
\end{eqnarray}
and
\begin{eqnarray}
r(k_x,k_y) =\frac{5\beta^2-1+(1+\beta^2)\cos(k_x\Delta x)}{5(1+\beta^2)-(5-\beta^2)\cos(k_x\Delta x)}\cos(k_y\Delta y)
\label{e:opt:lsor:hoc:bc}
\end{eqnarray}
and in case of any positive root for equation \refp{e:bcRobin:hyper:kx}, all the $\cos(k_x\Delta x)$ terms must be replaced by $\cosh(k_x\Delta x)$.

The required changes for the point SOR method of the HOC scheme in \refp{e:wopt:hoc2} and \refp{e:k2} are
\begin{align}
p &= \cos(k_x\Delta x)\;, \quad q = \cos(k_y\Delta y) \nonumber\\
k_1 &= \delta(k_x,k_y)\times
\left\{
\begin{array}{lr}
	\dfrac{1+10c_1}{\sqrt{84c_1^2+20c_1+1}}
	& c_1\geq 0 \\	
	1
	& c_1<0
\end{array}\right. \nonumber\\
k_2 &= -\frac{k_1^2}{2A}\times\frac{R_k A^2+A k_1+2R_k D}{R_k + k_1} \;,\quad	
R_{k} = \sqrt{k_1^2-\delta(k_x,k_y)^2}
\label{e:hocbc}
\end{align}
where in case of any positive root for equation \refp{e:bcRobin:hyper:kx}, in addition to replacing $\cos(k_x\Delta x)$ by $\cosh(k_x\Delta x)$, $k_x^2$ must be replaced by $-k_x^2$ in $\delta(k_x,k_y)$. 

\begin{figure}[htbp]
\setkeys{Gin}{draft=\showeps}
\includegraphics[width=0.45\textwidth]{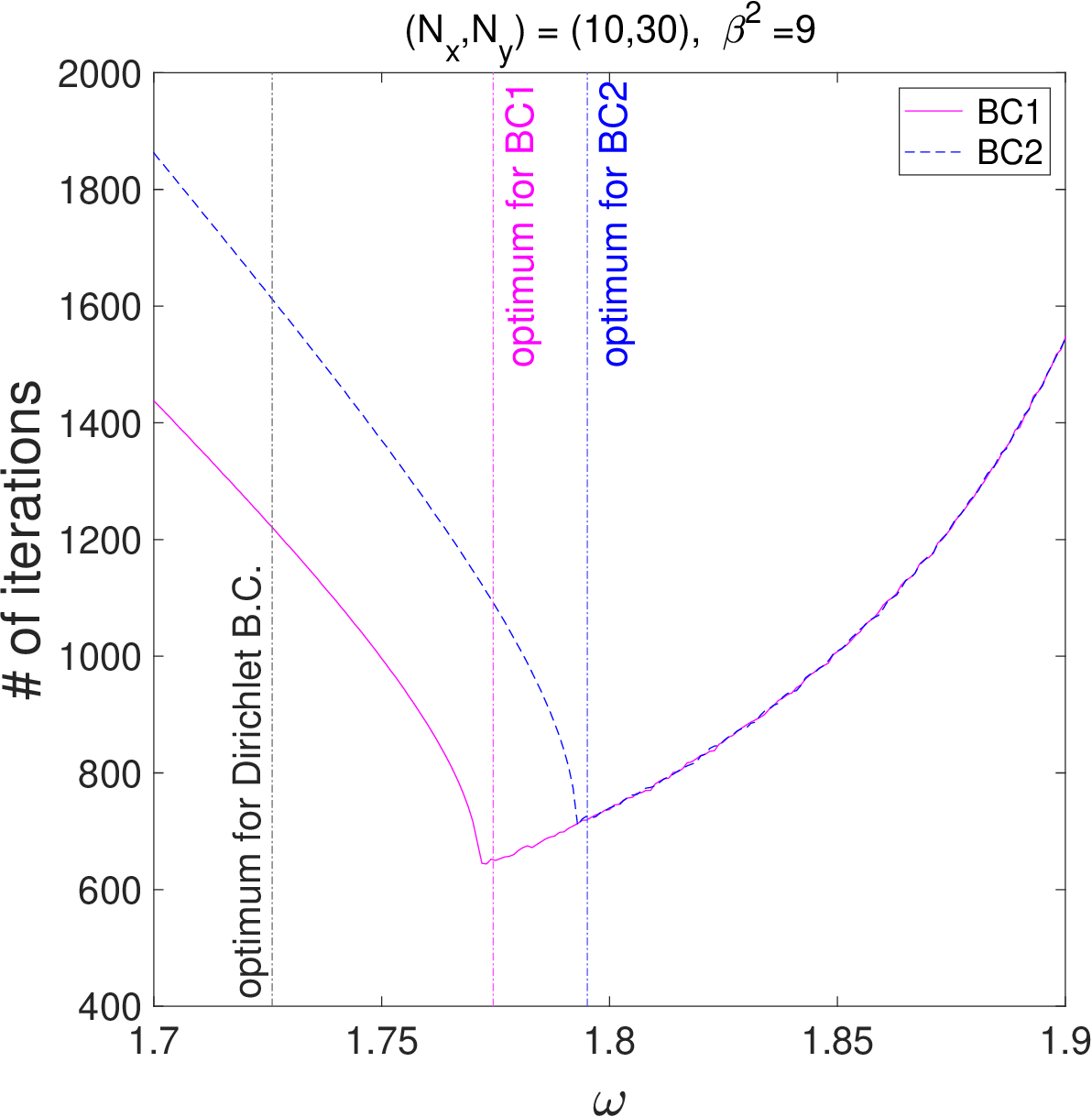}
\includegraphics[width=0.45\textwidth]{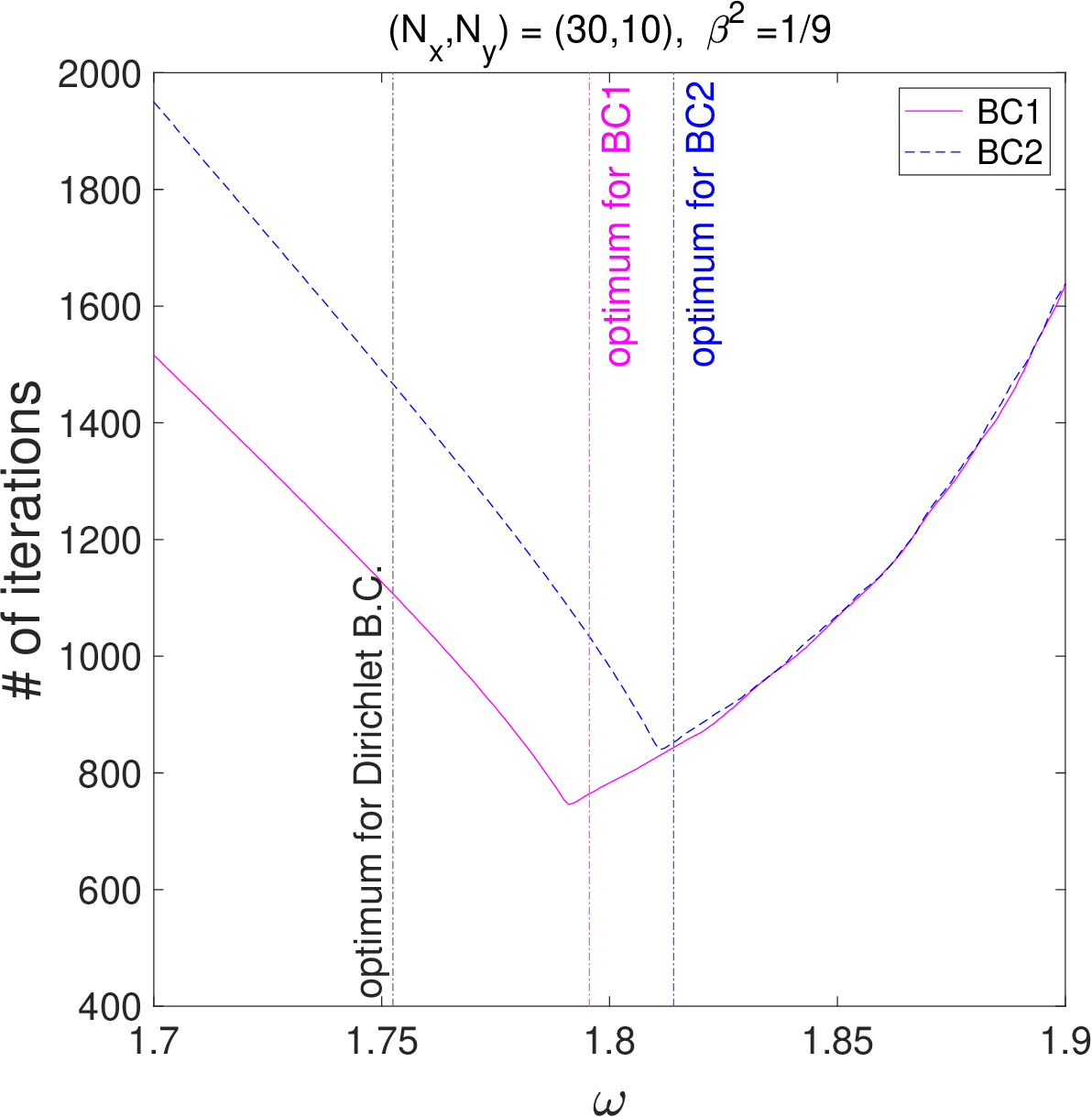}
\caption{Number of iterations for the point SOR method of the HOC scheme \refp{e:compact:sor} for different values of $\omega$. BC1: Neumann B.C. for the right edge. BC2: Neumann B.C. for both the left and right edges. The dash-dotted line indicates the optimal values for each case.}
\label{fig:plot:iter:hoc:bc}
\end{figure}

\begin{figure}[htbp]
\setkeys{Gin}{draft=\showeps}
\includegraphics[width=0.45\textwidth]{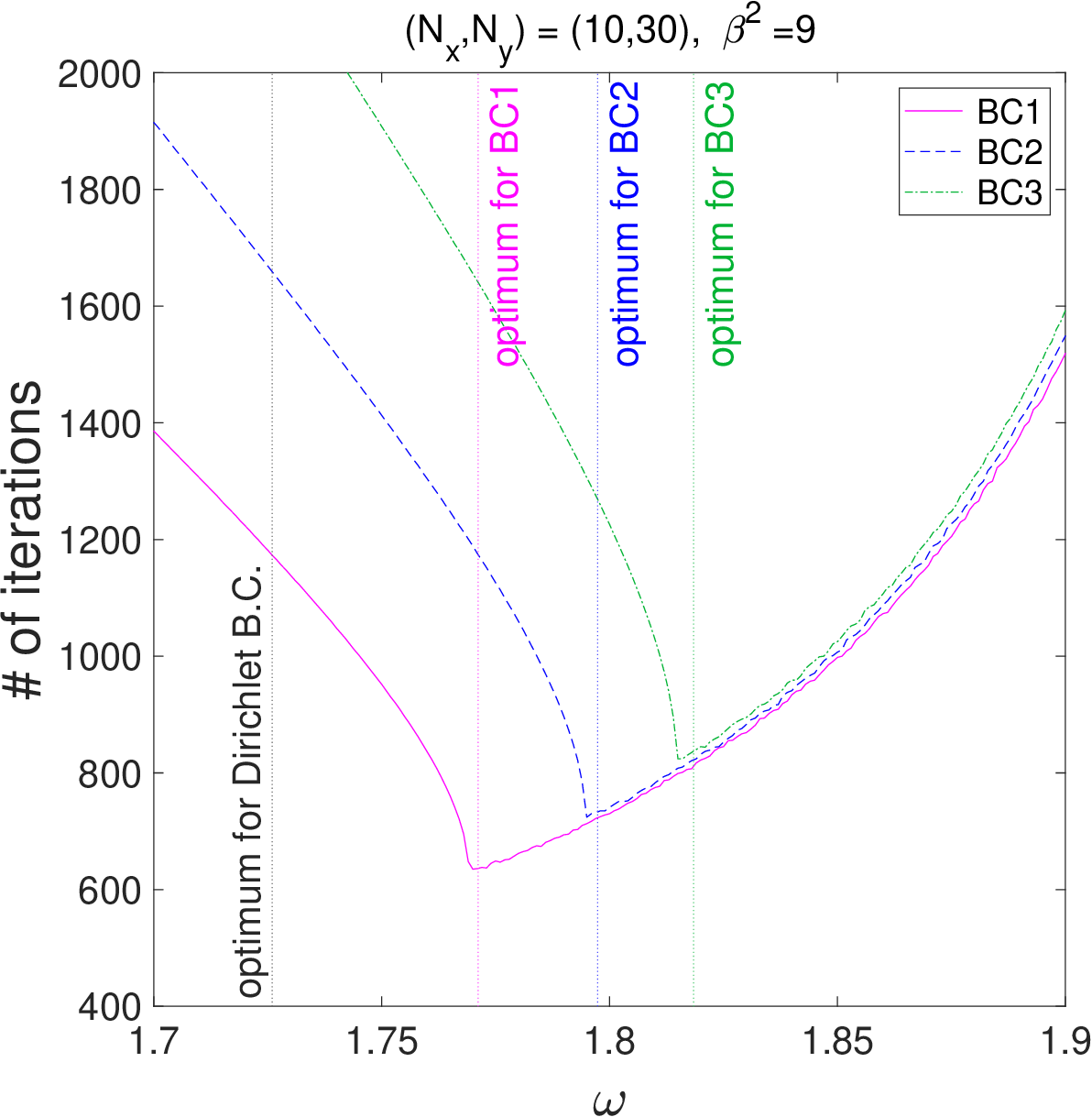}
\includegraphics[width=0.45\textwidth]{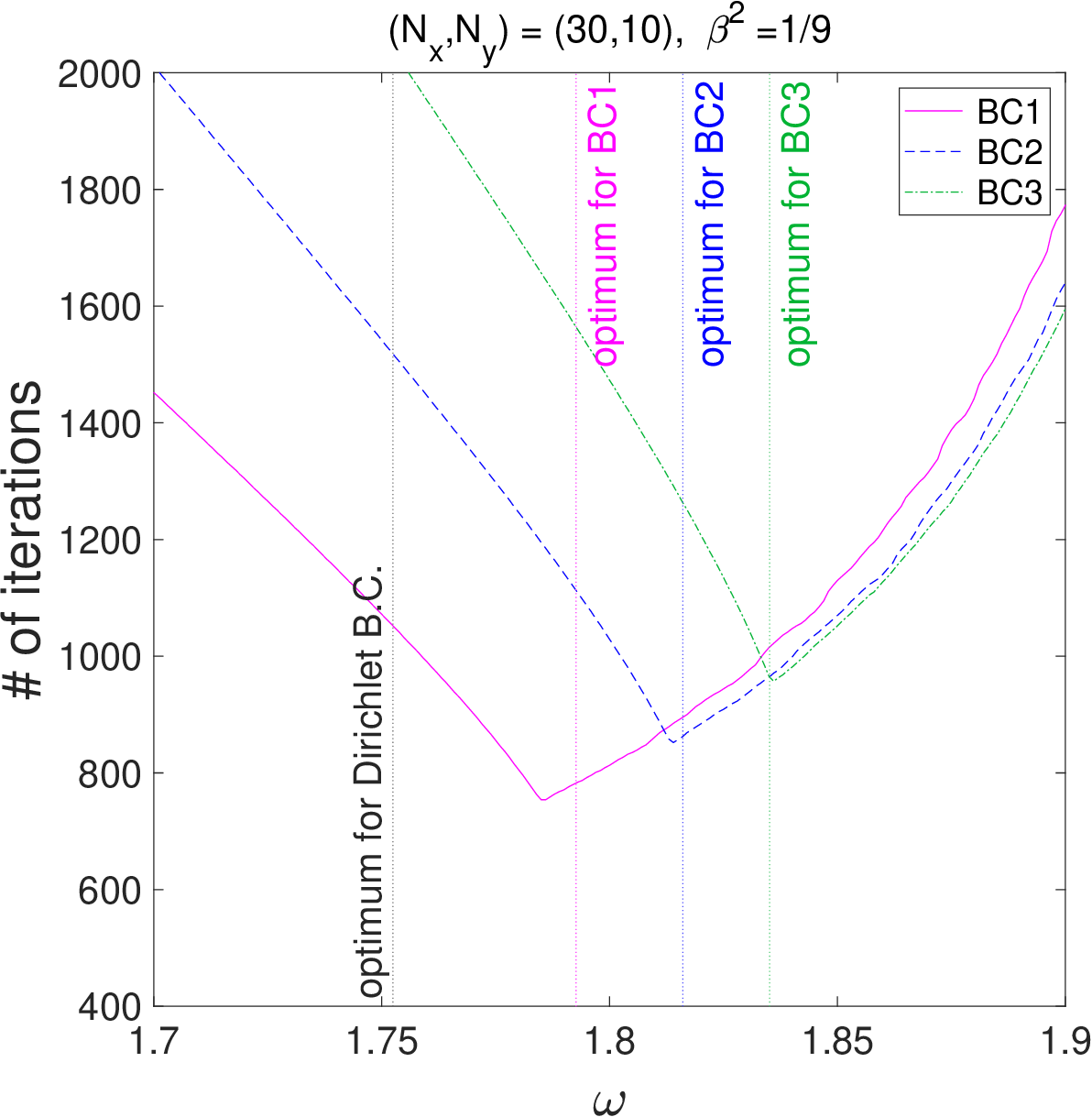}
\caption{Number of iterations for the point SOR method of the HOC scheme \refp{e:compact:sor} for different cases of the Robin boundary conditions (Table \ref{t:BCs}).}
\label{fig:plot:iter:hoc:robin}
\end{figure}

Figure \ref{fig:plot:iter:hoc:bc} shows the number of iterations required for the HOC scheme. The results verify the optimal value obtained by \refp{e:hocbc}. Again, we observe the number of iterations using the optimal value for Dirichlet boundary conditions \refp{e:wopt:hoc2}-\refp{e:k2} is nearly twice more than that of the actual optimal relaxation parameter. Figure \ref{fig:plot:iter:hoc:robin} shows the number of iterations for several cases of the Robin boundary conditions where the coefficients are given in Table \ref{t:BCs}.

\section{Conclusions}
\label{conclusions}
In the present study, we obtained the optimal relaxation parameter for solving the Poisson's equation with different types of boundary conditions on a rectangular grid using the SOR method. The emphasis was on grids with unequal grid step sizes in $x$- and $y$- directions and on Neumann and Robin boundary conditions. The central second-order and high-order compact schemes considered for the discretization. In addition to the point SOR method, the optimal relaxation parameter was presented for the line SOR method. The largest eigenvalue of each scheme was obtained by using the separation of variables. Minimizing the largest eigenvalue over the relaxation parameter led to a second-order equation except for the point SOR method of the HOC scheme which resulted in a quartic equation. For the latter, a first-order perturbation analysis was performed to find the optimal relaxation parameter. Due to low accuracy, the second-order terms were added to the perturbation analysis. Finally, all the obtained optimal relaxation parameters were verified numerically. 
Table \ref{t:summary} summarizes the results.
\begin{table}[h]
\footnotesize
\begin{center}
	\caption{Optimal relaxation parameters ($k_x$ and $k_y$ are given in Table \ref{t:2ndbcsummary}).}
	\label{t:summary}
	\begin{tabular}{|c|c|c|c|}
		\hline\hline
		\multicolumn{4}{|c|}{Second-order scheme} \\ \hline
		Method &        $\omega_{opt}$        &                            note                             &   largest eigenvalue   \\ \hline
		Point SOR  & $\frac{2}{1 + \sqrt{1-r^2}}$ &  $r=\frac{\CX+\beta^2\CY}{1+\beta^2}$  & $\omega_{opt}-1$ \\ \hline
		Line SOR  & $\frac{2}{1 + \sqrt{1-r^2}}$ & $r=\frac{\beta^2}{1+\beta^2 - \CX}\CY$ & $\omega_{opt}-1$ \\ \hline
		\hline
		\multicolumn{4}{|c|}{High-order compact scheme} \\ \hline
		Method& $\omega_{opt}$ &	note & largest eigenvalue \\\hline
		Point SOR   & $2-k_1h-k_2h^2$ & $k_1$ and $k_2$ are given by \refp{e:hocbc} & $1-2\beta_m h-2\gamma_m h^2$ \\\hline
		Line SOR  & $\frac{2}{1 + \sqrt{1-r^2}}$ & $	r =\frac{5\beta^2-1+(1+\beta^2)\CX}{5(1+\beta^2)-(5-\beta^2)\CX}\CY$ & $\omega_{opt}-1$ \\\hline		
	\end{tabular}
\end{center}
\end{table}

\section*{Acknowledgments}
The authors would like to thank the University of Tehran for financial support for this research under grant number 03/1/28745.

\appendix
\section{Number of roots of equation \refp{e:bcRobin:hyper:kx}}
\label{s:app}
Here, we investigate the existence of positive roots for equation \refp{e:bcRobin:hyper:kx}. Therefore, in the following we assume $k_x\geq 0$. First, by dividing the equation by $(ad-bc)\cosh k_x$ and dropping the subscript $x$, it is rewritten as  
\begin{eqnarray}
	\label{e:nmh}
	f(k):=(m-n \, \shk^2) \thk+ \shk =0 
\end{eqnarray}
where
\begin{eqnarray*}
	\shk := \frac{\sinh (k\Delta x)}{\Delta x} \;, \quad \thk := \tanh(k)
\end{eqnarray*}

The function $f$ is continuous and also an odd function. Consequently, its all even derivatives at $k=0$ are zero. We require the first and second derivatives for determining the number of roots:
\begin{align}
	\label{e:nmhd}
	f'(k)&=\shfk (1-2n \, \shk \thk)+(m-n \, \shk^2) \, \thfk  \\
	\label{e:nmhd2}
	f''(k)&=\shsk (1-2n \, \shk \thk) - 4n\, \shk\shfk \, \thfk \nonumber\\
	&-2n\,\shfk^2\,\thk - 2(m-n \, \shk^2)\thk \, \thfk
\end{align}
where
\begin{eqnarray*}
	\shfk = \cosh (k\Delta x)\;,\quad \shsk = \Delta x \sinh (k\Delta x) \;, \quad \thfk = 1-\thk^2
\end{eqnarray*}

Now, since $f$ is continuous and  
\begin{eqnarray}
	\label{e:nmh2}
	&f(0) = 0 \;, \qquad f'(0) = m+1 \\
	&\lim\limits_{k \rightarrow +\infty}f(k) = \text{sign}(n)\times -\infty
\end{eqnarray}
we conclude for $n(m+1)>0$, using the intermediate value theorem, $f$ has at least one positive root. 

In the following, to determine the sign of $f$ or its derivatives, we use the following inequalities
\begin{eqnarray*}
	\shk \geq 0\;,&\ \quad
	\cosh(k\Delta x) \geq 1 \\ 
	0\leq\thk \leq 1\;,& \quad
	0\leq \thfk \leq 1
\end{eqnarray*}
and also
\begin{eqnarray*}
	\shk \leq \sinh k \;, \quad \mathrm{for} \quad \Delta x \leq 1
\end{eqnarray*}
\subsection*{Case I: $n\leq 0$ and $m+1 > 0$}
Since in this case $n\leq 0$, we write the first derivative as below
\begin{align}
	\label{e:nmh3}
	f'(k) &= \cosh(k\Delta x) (1+2|n| \, \shk \, \thk)+(m+|n| \, \shk^2) \, \thfk \nonumber \\
	& \geq  \cosh(k\Delta x) + m \, \thfk \geq 1 + \min(m,0)
\end{align}
thus if $m+1> 0$, then the first derivative is always positive and the function is increasing which means it has no positive root. 

\subsection*{Case II: $n\leq 0$ and $m+1< 0$}
Since both $n$ and $m$ are assumed negative in this case, we have
\begin{align}
	\label{e:nmh4}
	f''(k)&=\shsk (1+2|n| \, \shk \thk) + 4|n|\, \shk\shfk \, \thfk \nonumber\\
	&+2|n|\shfk^2\thk - 2[-|m|+|n| \, \shk^2]\,\thk \, \thfk \nonumber\\
	&>2|n|\shfk^2\thk - 2[-|m|+|n| \, \shk^2]\,\thk \, \thfk \nonumber\\
	&>2|n|\shfk^2\thk - 2|n| \, \shk^2\,\thk \, \thfk \nonumber\\
	&=2|n|\,\thk \left[\shfk^2 - \shk^2 \, \thfk\right] \geq 0 
\end{align}
where in the last line, the term in the brackets can be shown that is positive. Therefore, the second derivative is positive and $f$ intersects $x$-axis exactly once.

\subsection*{Case III: $n > 0$ and $m+1>0$}
Since $m+1>0$, then $f$ is increasing at $k=0$ and $f>0$ immediate after $k=0$. However, for large $k$ the function becomes negative and there is at least one root. Now, we show if $f\leq 0$ then $f'<0$ which means after the root, the function is always decreasing and no other root exists.
If we assume $f'(k)<0$, then using the function definition \refp{e:nmh}, we obtain
\begin{eqnarray}
	\label{e:nmh5}	
	- n \leq -\frac{m}{\shk^2} - \frac{1}{\shk\,\thk}
\end{eqnarray}
now using this inequality for $f'$ in \refp{e:nmhd}, gives
\begin{eqnarray}
	\label{e:nmh6}	
	f'(k) \leq - 2m\,\shfk\frac{\thk}{\shk} - \shfk - \shk\frac{\thfk}{\thk}  
\end{eqnarray}
and since $-m<1$, we have
\begin{eqnarray}
	\label{e:nmh7}	
	f'(k) < 2\,\shfk\frac{\thk}{\shk} - \shfk - \shk\frac{\thfk}{\thk}  \leq 0 
\end{eqnarray}

\subsection*{Case IV: $n > 0$ and $m+1<0$}
If $k^*$ is a positive root of $f$, then  
\begin{eqnarray}
	\label{e:nmh8}	
	- n = -\frac{m}{\sshk^2} - \frac{1}{\sshk\,\sthk}
\end{eqnarray}
and replacing \refp{e:nmh8} in \refp{e:nmhd}, gives
\begin{align*}
	f'(k^*) &= - 2m\,\sshfk\frac{\sthk}{\sshk} - \sshfk - \sshk\frac{\sthfk}{\sthk}  \nonumber\\
	&= - 2(m+1)g(k^*)+ h(k^*)
\end{align*}
where
\begin{align}
	\label{e:nmh9}	
	g(k)&:= \shfk\frac{\thk}{\shk}\nonumber\\
	h(k)&:= 2g(k) - \shfk - \shk\frac{\thfk}{\thk}
\end{align}

The functions $g$ and $h$ are decreasing functions. Now, since $m+1<0$, then we have
\begin{eqnarray*}
	f'(k^*) = 2|m+1|g(k^*) + h(k^*)
\end{eqnarray*}
now, if $k_1^*$ and $k_2^*$ are two distinct roots and $k_1^*<k_2^*$, then
\begin{eqnarray*}
	f'(k_2^*) < f'(k_1^*)
\end{eqnarray*}

This inequality means there are two roots at most. Because, the continuity of $f$ means the sign of the function slope at roots alternately changes.

\subsection*{Case V: $m+1=0$}
In this case, we have
\begin{eqnarray*}
	f(0) = f'(0)= f''(0) = 0 \;, \quad f'''(0) = 2-6n+\Delta x ^2
\end{eqnarray*}
which means $k=0$ is a third-order root. Now, if $f'''(0)>0$ or equivalently $n<\frac{1}{3}+\frac{\Delta x ^2}{6}$, then for $n\leq 0$ this case is similar to case I and for $0<n<\frac{1}{3}+\frac{\Delta x ^2}{6}$ it is similar to case III. If $f'''(0)<0$, it is similar to case IV. For $f'''(0)=0$, we have $f''''(0)=0$ and $f'''''(0)<0$, and therefore the function is decreasing at $k=0$. Also, using \refp{e:nmh6} shows the function is decreasing for $k>0$ and thus it has no root. Note that, since we are interested in $1+4mn \geq 0$, therefore, $n\leq \frac{1}{4}$ and we always have $f'''(0)>0$ for this case.

\bibliographystyle{unsrt} 
\bibliography{mybibfile}

\begin{thebibliography}{10}

\bibitem{Young1971-va}
David~M Young.
\newblock {\em Iterative solution of large linear systems}.
\newblock Computer Science and Applied Mathematics. Academic Press, San Diego,
  CA, September 1971.

\bibitem{saad2003iterative}
Y.~Saad.
\newblock {\em Iterative Methods for Sparse Linear Systems}.
\newblock EngineeringPro collection. Society for Industrial and Applied
  Mathematics, 2003.

\bibitem{Liu2021}
Chein-Shan Liu.
\newblock A feasible approach to determine the optimal relaxation parameters in
  each iteration for the sor method.
\newblock {\em Journal of Mathematics Research}, 13:1, 01 2021.

\bibitem{Bai2003}
Zhong zhi Bai and Xue bin Chi.
\newblock Asymptotically optimal successive overrelaxation methods for systems
  of linear equations.
\newblock {\em Journal of Computational Mathematics}, 21(5):603--612, 2003.

\bibitem{MENG2014707}
Guo-Yan Meng.
\newblock A practical asymptotical optimal sor method.
\newblock {\em Applied Mathematics and Computation}, 242:707--715, 2014.

\bibitem{leveque1988sor}
Randall~J. Leveque and Lloyd~N. Trefethen.
\newblock Fourier analysis of the sor iteration.
\newblock {\em IMA Journal of Numerical Analysis}, 8(3):273--279, 07 1988.

\bibitem{adams1988sor}
Loyce~M. Adams, Randall~J. Leveque, and David~M. Young.
\newblock Analysis of the sor iteration for the 9-point laplacian.
\newblock {\em SIAM Journal on Numerical Analysis}, 25(5):1156--1180, 1988.

\bibitem{Tony1989}
Tony~F. Chan and Howard~C. Elman.
\newblock Fourier analysis of iterative methods for elliptic problems.
\newblock {\em SIAM Review}, 31(1):20--49, 1989.

\bibitem{kuo1990two}
C-C~Jay Kuo and Tony~F Chan.
\newblock Two-color fourier analysis of iterative algorithms for elliptic
  problems with red/black ordering.
\newblock {\em SIAM Journal on Scientific and Statistical Computing},
  11(4):767--793, 1990.

\bibitem{YANG2009325}
Shiming Yang and Matthias~K. Gobbert.
\newblock The optimal relaxation parameter for the sor method applied to the
  poisson equation in any space dimensions.
\newblock {\em Applied Mathematics Letters}, 22(3):325--331, 2009.

\bibitem{spotz1995high}
William~Frederick Spotz.
\newblock {\em High-order compact finite difference schemes for computational
  mechanics}.
\newblock The University of Texas at Austin, 1995.

\bibitem{hoffmann2000computational}
K.A. Hoffmann and S.T. Chiang.
\newblock {\em Computational fluid dynamics}, volume~2.
\newblock Engineering Education System, 4 edition, 2000.

\bibitem{Axelsson_1994}
Owe Axelsson.
\newblock {\em Iterative Solution Methods}.
\newblock Cambridge University Press, 1994.

\end{thebibliography}

\end{document}